\documentclass[12pt,a4paper]{amsart}
\usepackage{amsmath}
\usepackage{amssymb}
\usepackage{amsthm}
\usepackage{latexsym}
\usepackage{pstricks}
\usepackage{amsbsy}

\usepackage{array}
\usepackage[all]{xy}
\usepackage{amscd}
\usepackage{mathrsfs}
\usepackage{txfonts}
\usepackage{amsfonts}
\usepackage{graphicx}
\usepackage{listings}
\usepackage{url}
\usepackage{bookmark}
\newtheorem{theorem}{Theorem}[section]
\newtheorem{lemma}[theorem]{Lemma}
\newtheorem{corollary}[theorem]{Corollary}
\newtheorem{proposition}[theorem]{Proposition}
\newtheorem{definition}[theorem]{Definition}
\newtheorem{remark}[theorem]{Remark}

\newcommand{\abs}[1]{\left\lvert#1\right\rvert}

\definecolor{darkgreen}{rgb}{0.33, 0.42, 0.18}

\providecommand{\Gen}{\mathop{\rm Gen}\nolimits}%
\providecommand{\cone}{\mathop{\rm cone}\nolimits}%
\providecommand{\rad}{\mathop{\rm rad}\nolimits}%
\providecommand{\coker}{\mathop{\rm coker}\nolimits}%

\def\C{\mathcal{C}}

\def\P{\mathcal{P}}
\def\X{\mathbb{X}}
\def\Y{\mathbb{Y}}
\def\E{\mathcal{E}}

\def\K{\mathcal{K}}
\def\F{\mathcal{F}}
\def\T{\mathcal{T}}
\def\U{\mathcal{U}}

\def\W{\mathsf{W}}
\def\PP{{\mathbb P}}
\providecommand{\add}{\mathop{\rm add}\nolimits}%
\providecommand{\End}{\mathop{\rm End}\nolimits}%
\providecommand{\Ext}{\mathop{\rm Ext}\nolimits}%
\providecommand{\Hom}{\mathop{\rm Hom}\nolimits}%
\providecommand{\ind}{\mathop{\rm ind}\nolimits}%
\newcommand{\module}{\mathop{\rm mod}\nolimits}%

\usepackage{marginnote}
\newcommand{\sbt}{\,\begin{picture}(-1,1)(-1,-3)\circle*{3}\end{picture}\ }
\begin{document}

\title{$\tau$-exceptional sequences}

\author[Buan]{Aslak Bakke Buan}
\address{
Department of Mathematical Sciences \\
Norwegian University of Science and Technology \\
7491 Trondheim \\
Norway \\
}
\email{aslak.buan@ntnu.no}

\author[Marsh]{Bethany Rose Marsh}
\address{School of Mathematics \\ 
University of Leeds \\ 
Leeds, LS2 9JT \\ 
United Kingdom \\
}
\email{B.R.Marsh@leeds.ac.uk}

\dedicatory{Dedicated to the memory of Ragnar-Olaf Buchweitz}


\begin{abstract}
We introduce the notions of $\tau$-exceptional and signed $\tau$-exceptional sequences for any finite dimensional algebra.
We prove that for a fixed algebra of rank $n$, and for 
any positive integer $t \leq n$, there is a bijection between the set of
signed $\tau$-exceptional sequences of length $t$, and (basic) ordered support $\tau$-rigid objects with $t$ indecomposable direct summands.
If the algebra is hereditary, our notions coincide with exceptional and signed exceptional sequences. The latter were recently introduced by Igusa and Todorov, who constructed a similar bijection in the hereditary setting.
\end{abstract}


\thanks{
This work was supported by FRINAT grant numbers 231000 and 301375 from the Norwegian Research Council, by the Engineering and Physical Sciences
Research Council [grant number EP/G007497/1], and by the Mittag-Leffler Institute
(\emph{Representation Theory} programme, 2015). Part of the work for this paper was done while B. R. Marsh was a Guest Professor at NTNU, Trondheim, August--December 2014.}

\maketitle
\section*{Introduction}
An object in an abelian or triangulated category is said to be \emph{exceptional} if its
endormorphism algebra is a division ring and it has no self-extensions of any degree. Exceptional sequences are sequences of exceptional objects satisfying certain orthogonality conditions involving the vanishing of $\Hom$ and $\Ext$-groups. They were first introduced in an algebraic geometry setting~\cite{bondal,gorodentsev,gr} (see also~\cite{rud}). This motivated their consideration in the context of the representation theory of finite dimensional hereditary algebras (such as path algebras of quivers)~\cite{cb,r}.

The length of an exceptional sequence of finite-dimensional modules 
over an algebra 
is at most the number of simple modules.
If we have equality, the sequence is said to be \emph{complete}. 

Complete sequences always exist over hereditary algebras. Moreover,
there is in this case a transitive braid group action~\cite{cb,r}, making it possible to 
find all exceptional sequences and hence also all exceptional modules. Hence one 
can hope for a classification of such modules using combinatorial methods.

However, for an arbitrary finite-dimensional algebra, complete sequences do not always exist. For example, this is the case for the path algebra of an oriented two-cycle, modulo the square of its radical: only the projective indecomposable modules are exceptional. Since there are non-trivial morphisms between
the projective modules, they cannot form an exceptional sequence, and hence
no complete exceptional sequence exists for this algebra.
There are however non-hereditary settings where exceptional sequences are of importance; see~\cite{krause, hilleploog}.

Our primary aim in this article is to give an alternative approach to generalizing the definition of exceptional sequences from the hereditary to the general case,
ensuring that {\em complete} sequences always exist for arbitrary finite dimensional algebras.
Our work is firstly motivated by the signed exceptional sequences for hereditary finite dimensional algebras $H$ which
were recently introduced by Igusa and Todorov~\cite{it}. In this case, the projective
objects appearing in the sequence can be signed. Such sequences were needed
in order for the authors to define the \emph{cluster morphism category} of $H$, whose objects are the finitely generated wide subcategories of $\module H$.
The article~\cite{it} shows that for a hereditary algebra of finite representation type, the geometric realization of the cluster morphism category is an Eilenberg-MacLane space whose fundamental group coincides with the picture group of~\cite{itw}.

Signed exceptional sequences were needed to explain the composition and
associativity of maps in the cluster morphism category.
In particular it was shown that complete signed exceptional sequences are
in bijection with ordered cluster-tilting objects in the cluster category~\cite{bmrrt} corresponding to $H$. These are known to be in bijection with ordered clusters in the corresponding (acyclic) cluster algebra~\cite{bmrt,ck}.

Secondly, Adachi, Iyama and Reiten~\cite{air}, motivated by the behaviour of
cluster categories and cluster-tilted algebras (or, more generally, $2$-Calabi-Yau
categories and $2$-Calabi-Yau-tilted algebras), recently introduced \emph{$\tau$-tilting theory} for finite-dimensional algebras (note that some related ideas also 
appear in~\cite{derksenfei}). This included, in particular, the notions of
$\tau$-rigid modules and support $\tau$-tilting objects.
This theory has been explored in a number of subsequent papers.

Motivated by these ideas, we introduce the notion of a (signed) $\tau$-exceptional sequence of modules over a finite dimensional algebra. Although a module $M$ appearing in such a sequence is not exceptional in general, it does satisfy $\Ext^1(M,M)=0$ (using the
Auslander-Reiten formula: see Corollary~\ref{c:moduleproperty}).
We show that any finite dimensional algebra has a complete (signed or unsigned)
$\tau$-exceptional sequence. Furthermore, the complete signed $\tau$-exceptional sequences are in bijection with ordered support $\tau$-tilting objects (see Corollary~\ref{main-cor}),
generalizing the result of Igusa and Todorov in the hereditary case. Our approach is very different from that in \cite{it}. As well as dealing with the general case, an important ingredient in our proof is the  correspondence~\cite{air} between $\tau$-rigid modules and rigid $2$-term 
complexes in the derived category.

In the case of a cluster-tilted algebra~\cite{bmr}, there is a connection to cluster algebras, which we now describe. Fix a cluster-tilted algebra, i.e.\ an algebra of the form $\Lambda=\End_{\C_Q}(T)$, where $Q$ is a finite acyclic quiver, $\C_Q$ is the corresponding cluster category, and $T$ is a cluster-tilting object in $\C_Q$.
By~\cite{air}, support $\tau$-tilting $\Lambda$-modules are in bijection with cluster-tilting objects in $\C_Q$.
By~\cite{bmrt,ck}, cluster-tilting objects in $\C_Q$ are in bijection with
clusters in the cluster algebra $\mathcal{A}_Q$~\cite{fz} associated to $Q$.
Both of these bijections induce bijections between the ordered versions.
It follows from these results and the results proved here that there is a bijection
between complete signed $\tau$-exceptional sequences over $\Lambda$ and ordered
clusters in the cluster algebra $\mathcal{A}_Q$.
We remark that if $\Lambda$ is hereditary, (signed) $\tau$-exceptional sequences and (signed) exceptional sequences coincide (although in general they are different), and we obtain
a bijection between complete signed exceptional sequences over $\Lambda$ and
ordered clusters in $\mathcal{A}_Q$; in this case such a bijection was also obtained
in~\cite{it}.

A finite dimensional algebra is called {\em $\tau$-tilting finite} \cite{dij} if there are only finitely many basic $\tau$-tilting objects (up to isomorphism).
Motivated by the cluster morphism categories mentioned above, in~\cite{bm} we construct a natural category whose objects are the wide subcategories of the module category of a $\tau$-tilting finite algebra. The construction relies heavily on the bijection established in this paper.
After the first preprint-version of this paper and \cite{bm} appeared,  Hanson and Igusa~\cite{hi} used our results to show that a generalized version of the picture group can be understood in terms of the category of wide subcategories
from~\cite{bm}.

The paper is organized as follows. In Section \ref{sec-not} we give some 
notation and background, and state the main result. In Section \ref{sec:two-term} we 
give some background and results concerning 2-term silting objects in
the derived category,
and their links to $\tau$-rigid objects in module categories. 
    In Sections \ref{sec-exc} and \ref{sec-red} we prepare
    for the proof of the main result, which is then completed in Section \ref{sec-main}. We conclude by giving some examples in Section \ref{sec-ex}.

Both authors would like to thank the referee for a thorough reading of the first version of the paper, and their helpful comments, and would like to thank William Crawley-Boevey, Kiyoshi Igusa, Osamu Iyama, Gustavo Jasso, Hugh Thomas and Gordana Todorov for stimulating discussions. In particular, the authors want to thank Iyama and Jasso for sharing the idea which led to Lemma \ref{exchange}. 
ABB would like to thank BRM and the School of Mathematics at the University of Leeds for their warm hospitality, and BRM would like to thank ABB and the Department of Mathematical Sciences at NTNU for their warm hospitality.

\section{Notation and main result}\label{sec-not}
Let $\Lambda$ be a basic finite dimensional algebra over a field
$k$, and let $\module \Lambda$ denote the category of finite
dimensional $\Lambda$-modules. Let $\P(\Lambda)$ denote the full
subcategory of projective objects in $\module \Lambda$. 
Similarly, if $\mathcal{X}$ is a subcategory of $\module \Lambda$, let
$\P(\mathcal{X})$ denote the full subcategory of $\mathcal{X}$ consisting of the
Ext-projective objects in $\mathcal{X}$, i.e. the objects $P$ in $\mathcal{X}$ such that $\Ext^1(P,X)=0$ for all $X\in \mathcal{X}$.

For an additive category $\C$ and an object $X$ in $\C$, we denote by
$\add X$ the additive subcategory of $\C$ generated by $X$, i.e. the full subcategory of $\C$ whose objects are all direct summands of direct sums of
copies of $X$.
For a subcategory $\mathcal{X} \subseteq \mathcal{C}$, we define
$\mathcal{X}^{\perp} = \{Y \in \C \mid \Hom(X,Y)= 0 \text{ for all } X \in \mathcal{X} \}$,
and define $^{\perp}\mathcal{X}$ similarly.

If $\C$ is (skeletally small) and Krull-Schmidt, we denote by $\ind(\C)$ the set of isomorphism classes of indecomposable objects in $\C$. For an object $X$, we write $\ind X$ for $\ind(\add(X))$. For any basic object $X$ in $\C$, let $\delta(X)$ denote the number of indecomposable direct summands of $X$. We denote $\delta(\Lambda)$ by $n$.

If $\C$ is abelian and $X$ is an object of $\C$, we denote by $\Gen X$ the full subcategory of $\C$ consisting of all objects which are factors of objects in
$\add X$.

In general, all subcategories considered are assumed to be full and closed under isomorphisms. All objects are taken to be basic where possible and considered up to isomorphism.

Let $\tau$ denote the Auslander-Reiten translate in $\module \Lambda$.
We now recall notation and definitions of from \cite{air}. Note that our definitions are 
slightly different, but clearly equivalent to the corresponding definitions in 
\cite{air}.

An object  $U$  in $\module \Lambda$ is called \emph{$\tau$-rigid} if
$\Hom(U,\tau U) = 0$. Let $D^b(\module \Lambda)$ denote the bounded derived category of
$\module \Lambda$, with shift functor denoted by $[1]$.
We consider $\module \Lambda$ as a full subcategory of $D^b(\module \Lambda)$ by regarding a module as a stalk complex concentrated in degree $0$.
Consider the full subcategory $\C(\Lambda) = \module \Lambda \amalg 
\module \Lambda [1]$ of $D^b(\module \Lambda)$.

\begin{definition}
The object $M \amalg P[1]$ in $\C(\Lambda)$ is called \emph{support $\tau$-rigid} if 
\begin{itemize}
\item[(i)] The object $M$ lies in $\module \Lambda$ and satisfies $\Hom(M, \tau M)= 0$, and
\item[(ii)] The object $P$ lies in $\P(\Lambda)$ and satisfies $\Hom(P,M)= 0$.
\end{itemize}
The object $M \amalg P[1]$ is called a \emph{support $\tau$-tilting object} if $\delta(M \amalg P[1])= n$.
Moreover, $M$ is in this case called a \emph{support $\tau$-tilting module} or just
a $\tau$-tilting module if in addition $P=0$.

\end{definition}

We want to consider all possible orderings of such objects, in the following sense.

\begin{definition}
For a positive integer $t$, an ordered $t$-tuple of
indecomposable objects 
$(\T_1, \dots, \T_t)$ in $\C(\Lambda)$ is called an \emph{ordered 
support $\tau$-rigid object} if $\amalg_{i=1}^t \T_i$ is a  
basic support $\tau$-rigid object.  If, in addition, $t=n$, then 
$(\T_1, \dots, \T_t)$ is called an ordered support $\tau$-tilting object.
\end{definition}

Let $U$ be an indecomposable $\tau$-rigid $\Lambda$-module. 
Jasso \cite{jasso} considered the category $J(U) = U^{\perp} \cap {^{\perp}(\tau U)}$, and, in particular, proved that $J(U)$ is equivalent to the module category of a finite-dimensional algebra $\Lambda_U$ with $\delta(\Lambda_U) = \delta(\Lambda)-1$. By~\cite[Cor. 3.22]{bst},~\cite[Theorem 3.28]{dirrt} $J(U)$ is an
exact abelian (wide) subcategory of $\module \Lambda$ (see Proposition \ref{jass} for more details).
We say that an object $N$ in $J(U)$ is \emph{$\tau$-rigid in $J(U)$} if the corresponding $\Lambda_U$-module is $\tau$-rigid.

Note that a $\Lambda$-module which is not $\tau$-rigid can be $\tau$-rigid in $J(U)$. For example, let $Q$ be the quiver 
$$\xymatrix@=5mm{
& 2  \ar^{\beta}[dr] & \\
1  \ar^{\alpha}[ur]  \ar_{\gamma}[rr]& & 3
} 
$$
and let $\Lambda = kQ/I$ where $I$ is the ideal generated by the path $\beta \alpha$.
Let $S_i$ (respectively, $I_i$, $P_i$) be the simple (respectively, indecomposable injective, indecomposable projective) module associated to vertex $i$ of $Q$. Then $\tau I_3=P_1$ and $I_3$ is not $\tau$-rigid in $\module\Lambda$.
However, $I_3$ is projective in $J(S_2)$ (which is equivalent to the module category of a quiver of Dynkin type $A_2$), so $I_3$ is $\tau$-rigid in $J(S_2)$. See Section~\ref{ex3} for more details.

For a projective object $P$, we let $J(P[1]) = J(P) = P^{\perp}$ and set $\Lambda_{P[1]}=\Lambda_P$.
Then, for any indecomposable object $\U$ in $\C(\Lambda)$, we have an equivalence $F_{\U}:J(\U)\rightarrow \module \Lambda_{\U}$.

For a full subcategory $\mathcal{Y}$ of $\module \Lambda$, we denote by $\C(\mathcal{Y})$ the full subcategory $\mathcal{Y} \amalg \mathcal{Y}[1]$ of $\C(\Lambda)$. Thus, for $\U$ as above, we have the full subcategory $\C(J(\U))$ of $\C(\Lambda)$. Since $J(\U)$ is a wide subcategory of $\module \Lambda$, we have that for $A,B\in J(U)$,
$$\Hom_{D^b(\Lambda)}(A,B[1])\cong \Ext^1_{\Lambda}(A,B)\cong
\Ext^1_{J(\U)}(A,B)\cong \Hom_{D^b(\Lambda_U)}(F_{\U}(A),F_{\U}(B)[1]).$$
It follows that $\C(J(\U))\simeq \C(\module \Lambda_{\U})$.

We say that an object $N\amalg Q[1]$ in $\C(J(\U))$ is \emph{support $\tau$-rigid in $\C(J(\U))$} if the corresponding object in $\C(\module \Lambda_{\U})$ is support $\tau$-rigid, i.e.\ $N$ is $\tau$-rigid in $J(\U)$ and $Q$ is an object in $J(\U)$ which is projective as an object in $J(\U)$.

This allows us to define signed $\tau$-exceptional sequences recursively as follows:

\begin{definition}
For a positive integer $t$, an ordered $t$-tuple of indecomposable objects 
$(\U_1, \dots, \U_{t-1}, \U_t)$ in $\C(\Lambda)$ is called a
\emph{signed $\tau$-exceptional sequence}, if 
$\U_t$ is a support $\tau$-rigid object in $\C(\Lambda)$ and
$(\U_1, \dots, \U_{t-1})$ is a signed $\tau$-exceptional sequence in $\C(J(\U_t))$.
A signed $\tau$-exceptional sequence $(\U_1, \dots, \U_{t-1}, \U_t)$
in which all of the $\U_i$ lie in $\module \Lambda$ is called a \emph{$\tau$-exceptional sequence}. A (signed) $\tau$-exceptional sequence of length $n$ is said to be \emph{complete}.
\end{definition}

In the above example, $(P_2,I_3,S_2)$ is a $\tau$-exceptional sequence in $\module\Lambda$. Note that, although $I_3$ is not projective in $\module \Lambda$, it is projective in $J(S_2)$ and therefore can take a `sign'; in other words,
$(P_2,I_3[1],S_2)$ is a signed $\tau$-exceptional sequence in $\module \Lambda$
(see Section~\ref{ex3}). This is the reason we take $\module \Lambda[1]$ as a summand in the definition of $C(\Lambda)$ above, rather than just the shift of the subcategory of projective $\Lambda$ modules.

We can now state our main result.

\begin{theorem}\label{main-int}
Let $\Lambda$ be  a finite dimensional algebra. For each $t \in \{1, \dots, n \}$
there is a bijection between the set of ordered support $\tau$-rigid
objects of length $t$ in $\C(\Lambda)$ and the
set of signed $\tau$-exceptional sequences of length $t$ in $\C(\Lambda)$. 
\end{theorem}

For $t=n$, we obtain the following.

\begin{corollary}\label{main-cor-int}
Let $\Lambda$ be  a finite dimensional algebra. 
Then there is a bijection between the set of ordered support $\tau$-tilting
objects in $\C(\Lambda)$ and the
set of complete signed $\tau$-exceptional sequences in $\C(\Lambda)$. 
\end{corollary}

We recall the following crucial result of \cite{air} which shows that, given a $\tau$-rigid module, there are two support $\tau$-tilting objects containing the module as a direct summand, satisfying certain properties.

\begin{theorem}\cite[Section 2.3]{air}\label{two-pairs}
Let $U$ be a $\tau$-rigid module. 

\begin{itemize}
\item[(a)] Up to isomorphism, there is a unique basic module $B[U]$ such that $B[U] \amalg U$ is a $\tau$-tilting  module and $\add (B[U] \amalg U) = \P({^{\perp}{\tau U}})$. Moreover, we
have ${^{\perp}{\tau U}} = {^{\perp}{\tau (U \amalg B[U])}}  = \Gen(U \amalg B[U])$. The equality ${^{\perp}{\tau U}} = {^{\perp}{\tau (U \amalg B[U])}}$ characterizes $B[U]$.
\item[(b)]  Up to isomorphism, there is a unique basic module $C[U]$ and a basic projective module $Q$ such that $C[U] \amalg U \amalg Q[1]$ is a 
support $\tau$-tilting object and $\add (C[U] \amalg U)  = \P(\Gen U)$.  
In particular, we have
$\add Q = \P(\Lambda) \cap {^{\perp}{U}}$.
\end{itemize}
\end{theorem}

\begin{proof} 
Part (a), apart from the characterization, follows from~\cite[Theorem 2.10]
{air}.
For the characterization, suppose that $B$ satisfies ${^{\perp}{\tau U}} = {^{\perp}{\tau (U \amalg B)}}$ and that $U\amalg B$ is a basic $\tau$-tilting module. By~\cite[Theorem 2.12(c)]{air}, we have ${^{\perp}{\tau (U \amalg B)}}=\Gen
U \amalg B$.  By~\cite[Lemma 2.11]{air}, $^{\perp}{\tau U}$ is functorially finite, so
$U\amalg B$ is Ext-projective in $^\perp{\tau U}$ by~\cite[Prop 2.9]{air}, and hence
$B \amalg U$ is in $\add B[U] \amalg U$ by assumption. Now, since $B \amalg U$ and
$B[U] \amalg U$ are both basic $\tau$-tilting modules, we must have that
$B$ and $B[U]$ are isomorphic.

For part (b), let $C[U]$ be defined so that $C[U]\amalg U$ is basic and $\add(C[U]\amalg U)=\P(\Gen U)$.
Let $Q$ be basic and such that 
$\add Q=\P(\Lambda)\cap {}^{\perp} U$. Since $\Gen U$ is functorially finite, it follows from~\cite[Theorem 2.7]{air} that $C[U]\amalg U$ is support $\tau$-tilting. We have that $\Gen U=\Gen C[U]\amalg U$, since $C[U]$ is in $\Gen U$. If $P$ is a projective module it is easy to see that $\Hom(P,M)=0$ for a module $M$ if and only if $\Hom(P,N)=0$ for all $N\in \Gen M$. Hence, we have:
$$\P(\Lambda)\cap {}^{\perp} U=\P(\Lambda)\cap {}^{\perp} (\Gen U)=
\P(\Lambda) \cap {}^{\perp} (C[U]\amalg U)$$
and it follows that $C[U]\amalg U\amalg Q[1]$ is a support $\tau$-tilting object.
\end{proof}

The module $B[U]$ in Theorem~\ref{two-pairs} is known as the \emph{Bongartz complement} of $U$, and we refer to $C[U]$ as the \emph{co-Bongartz complement} of $U$.
An important step towards our main theorem is the construction of explicit bijections between the indecomposable direct summands of $B[U]$ and the indecomposable direct summands of $C[U] \amalg Q[1]$.

\section{2-term rigid objects}\label{sec:two-term}

In this section we 
discuss 2-term silting objects in
the derived category
and their links to $\tau$-rigid objects in the module category. 

We denote by $H^n$ the functor from $D^b(\module \Lambda)$ which maps a complex to its $n$th homology.
We regard the bounded homotopy category of 
projectives $\K= K^b(\P(\Lambda))$ as a full subcategory of $D^b(\module 
\Lambda)$.
An object $\mathbb{U}$ in $\K$ is said to be \emph{rigid} if 
$\Hom(\mathbb{U}, \mathbb{U}[i]) = 0$ for all $i > 0$, and \emph{silting} if in addition
it generates $\K$ as a triangulated category (i.e. $\mathbb{U}$ is not contained
in any proper triangulated subcategory of $\K$ closed under direct summands).

An object of the form $$\cdots \to 0 \to P^{-1} \to P^0 \to 0 \to \cdots$$
in $\K$ is called a \emph{$2$-term object}. If $P^0=0$, we denote such an object by $P^{-1}[1]$.
A $2$-term object in $\K$ which is rigid (respectively, silting) is called \emph{$2$-term rigid} (respectively, \emph{$2$-term silting}).
 
For a module $U$, we let $\PP_U$ denote the minimal projective presentation of
$U$, considered as a $2$-term object in $\K$.
The following two lemmas are well-known.

\begin{lemma}\label{lift}
\begin{itemize}
\item[(a)] Let $\mathbb{X}$ and $\mathbb{Y}$ be 2-term objects in $\K$.
Then $H^0$ induces an epimorphism 
$\Hom_{\K}(\mathbb{X},\mathbb{Y}) \to \Hom(H^0(\mathbb{X}), H^0(\mathbb{Y}))$
with kernel consisting of the maps factoring through $\add \Lambda[1]=\add(\Lambda\rightarrow 0)$.
\item[(b)] Let $X,Y,Z$ be in $\module \Lambda$. If any map in $\K$ from $\PP_X $ to $\PP_Z$ factors through $\PP_Y$, then any map from $X$ to $Z$ factors through $Y$.
\end{itemize}
\end{lemma}

\begin{proof}
Part (a) is straightforward, and part (b) is a direct consequence of part (a).
\end{proof}

\begin{lemma}\label{indec}
Let $X$ be in $\module \Lambda$. Then $\PP_X$ is indecomposable in $\K$
if and only if $X$ is an indecomposable module.
\end{lemma}

\begin{proof}
Straightforward.
\end{proof}

We also need the following facts from \cite{air}.

\begin{lemma}\cite[Section 3]{air}\label{rigid-rigid}
Let $U,X$ be in $\module \Lambda$.
\begin{itemize}
\item[(a)] $\Hom(U, \tau X)= 0$ if and 
only if $\Hom_{\K}(\PP_X, \PP_U[1]) = 0$. In particular, the module $U$ is 
$\tau$-rigid if and only if $\PP_U$ is rigid.
\item[(b)] If $U \amalg P[1]$ is a support $\tau$-rigid object in $\C(\Lambda)$, then 
$\PP_U \amalg P[1]$ is
a $2$-term rigid object in $\K$, and it is $2$-term silting if and only if $U \amalg P[1]$ 
is support $\tau$-tilting.
\item[(c)] Let $\mathbb{U}$ be a $2$-term rigid object in $\K$. Write $\mathbb{U}$ in the form $\mathbb{P}\amalg Q[1]$, where $Q$ is a projective module and $\mathbb{P}$ has no direct summand of the form $Q'[1]$, where $Q'$ is a projective module. Then $H^0(\mathbb{P})\amalg Q[1]$ is support $\tau$-rigid in $C(\Lambda)$. If $\mathbb{U}$ is $2$-term silting, then $H^0(\mathbb{P})\amalg Q[1]$ is support $\tau$-tilting in $C(\Lambda)$.
\item[(d)] The constructions of (b) and (c) give a bijection between
$2$-term silting objects in $\K$ and support $\tau$-tilting objects in
$\C(\Lambda)$.
\end{itemize}
\end{lemma}

Note that by Lemma \ref{rigid-rigid} any $2$-term silting object in $\K$ is, up to homotopy, of the form $\PP_X \amalg Q[1]$, where $\PP_X$ is a minimal projective
presentation for some $\Lambda$-module $X$.
In particular, this implies the following.

\begin{lemma}\label{nokernel}
	Let $\mathbb{U}$ be a $2$-term silting object in $\K$. Then $H^0(\mathbb{U})$
	is $\tau$-tilting if and only if $\mathbb{U}$ has no direct summand of the form $P \to 0$
\end{lemma}
\begin{proof}
Since $\mathbb{U}$ is $2$-term silting, we have $\mathbb{U} \simeq \mathbb{P}_U \amalg P[1]$
for a support $\tau$-tilting object $U \amalg P[1]$ in $\C(\Lambda)$ with $H^0(\mathbb{U})= U$.
And we have $P = 0$ if and only if $U$ is  $\tau$-tilting.  
\end{proof}

Parts (b) to (e) in the following Lemma are essentially contained
in \cite[Section 2]{ai}, but we include proofs for convenience.
\begin{lemma}\label{exchange}
Let $\mathbb{U}$ and $\mathbb{X}$ be $2$-term objects in 
$\K$ satisfying
\begin{itemize}
\item[(i)] $\mathbb{U} \amalg \mathbb{X}$ is rigid;
\item[(ii)] $\add \mathbb{U} \cap \add \mathbb{X} = 0$.
\end{itemize}
Let 
\begin{equation}\label{x-y-triangle}
\mathbb{Y} \overset{\beta}{\rightarrow}  \mathbb{U'} \overset{\alpha}{\rightarrow} \mathbb{X} \to  
\end{equation}
be a triangle in which $\alpha$ is a minimal right $\add \mathbb{U}$-approximation. Then the following hold.
\begin{itemize}
\item[(a)] The object $\mathbb{Y}$ is a $2$-term object if and only if the induced map 
$H^0(\alpha) \colon H^0(\mathbb{U'}) \to H^0(\mathbb{X})$ is an epimorphism.
\item[(b)] The object $\mathbb{Y} \amalg \mathbb{U}$ is rigid. 
\item[(c)] The morphism $\mathbb{Y} \overset{\beta}{\rightarrow} \mathbb{U'}$ is a minimal left $\add \mathbb{U}$-approximation. 
\item[(d)] $\add \mathbb{U} \cap \add \mathbb{Y} = 0$.
\item[(e)] The object $\mathbb{Y} \amalg\mathbb{U}$ is silting if and only if
$\mathbb{X} \amalg \mathbb{U}$ is silting.
\item[(f)] If $\mathbb{X} \amalg \mathbb{U}$ is $2$-term silting and $H^0(\alpha)$ is an epimorphism, then $\mathbb{Y} \amalg \mathbb{U}$ is $2$-term silting.
\end{itemize}
\end{lemma}

\begin{proof}
Note that we can assume that any differential appearing in an object in $\K$ is radical, i.e. that its image is contained in the radical of its target.
Let $\cone(\alpha)$ be the mapping cone of $\alpha$, which has the form:
$$\cdots \to 0 \to U^{-1} \to U^0 \amalg X^{-1} \to X^0\to 0 \to \cdots$$
Then $\mathbb{Y} \simeq \cone(\alpha)[-1]$. 
Clearly $ \cone(\alpha)[-1]$ is $2$-term if and only if the map 
$U^0 \amalg X^{-1} \to X^0$ is an epimorphism, which is equivalent to 
$U^0 \to X^0$ being a (split) epimorphism, since the map $ X^{-1} \to X^0$ is
radical. This holds
if and only if the induced map $H^0(\alpha) \colon H^0(\mathbb{U'}) \to H^0(\mathbb{X})$ is an epimorphism, giving (a).

We now prove (b).
We apply $\Hom(\mathbb{U}, \ )$ to (\ref{x-y-triangle}), obtaining the exact sequences
$$\Hom(\mathbb{U}, \mathbb{U'}[i] ) \to \Hom(\mathbb{U}, \mathbb{X}[i]) \to \Hom(\mathbb{U}, \mathbb{Y}[i+1]) \to  \Hom(\mathbb{U}, \mathbb{U'}[i+1]),$$
for all $i$. Since $\alpha$ is a right $\add \mathbb{U}$-approximation, 
the map $\Hom(\mathbb{U}, \mathbb{U'} ) \to \Hom(\mathbb{U}, \mathbb{X})$ is an epimorphism. Hence $\Hom(\mathbb{U}, \mathbb{Y}[i]) = 0$ for all $i >0$.

Applying $\Hom(\mathbb{X}, \ )$ to (\ref{x-y-triangle}), and considering the exact sequence
$$\Hom(\mathbb{X}, \mathbb{X}[i-1]) \to \Hom(\mathbb{X}, \mathbb{Y}[i]) \to  \Hom(\mathbb{X}, \mathbb{U}[i])$$
gives
$\Hom(\mathbb{X}, \mathbb{Y}[i]) = 0$ for all
$i >1$.

Applying $\Hom(\ , \mathbb{Y})$ to (\ref{x-y-triangle})
gives an exact sequence
$$ \Hom(\mathbb{U'}, \mathbb{Y}[i] ) \to \Hom(\mathbb{Y}, \mathbb{Y}[i]) \to \Hom(\mathbb{X}, \mathbb{Y}[i+1]),$$
and hence $\Hom(\mathbb{Y}, \mathbb{Y}[i])=$ for all $i>0$.

Finally, applying $\Hom(\ , \mathbb{U})$ to (\ref{x-y-triangle})
gives an exact sequence
$$ \Hom(\mathbb{U'}, \mathbb{U}[i] ) \to \Hom(\mathbb{Y}, \mathbb{U}[i]) \to \Hom(\mathbb{X}, \mathbb{U}[i+1]),$$
and hence $\Hom(\mathbb{Y}, \mathbb{U}[i])=$ for all $i>0$.
This finishes the proof of (b).

For part (c) consider the exact sequence
$$  \Hom(\mathbb{U'}, \mathbb{U}) \to \Hom(\mathbb{Y}, \mathbb{U}) \to \Hom(\mathbb{X}, \mathbb{U}[1]),$$
which is part of the long exact sequence obtained by applying 
$\Hom(\ , \mathbb{U})$ to (\ref{x-y-triangle}).
Since the last term vanishes, the first map is an epimorphism, and
hence the map $\mathbb{Y} \overset{\beta}{\rightarrow} \mathbb{U'}$ is a left 
$\add \mathbb{U}$-approximation. It is minimal since
$\add \mathbb{U} \cap \add \mathbb{X} = 0$.

Part (d) follows directly from the minimality of $\alpha$ and part (c).

Part (e) follows from part (b) and the existence of
the triangle~\eqref{x-y-triangle}, and part (f) is a direct consequence of parts (a) and (e).
\end{proof}

\section{Exchange}\label{sec-exc}

Let $U$ be a $\tau$-rigid $\Lambda$-module. Recall from Theorem \ref{two-pairs} the notation $B[U]$ for the Bongartz complement of $U$ and
$C[U] \amalg Q[1]$ for the co-Bongartz complement of $U$ in $\C(\Lambda)$.
We will denote $C[U]$ by $C$ in the sequel.
The aim of this section is to give an explicit bijection between the indecomposable direct summands in these two complements of $U$.

\begin{remark}\label{rema}
Let $\mathbb{C}_Q = \PP_{C} \amalg Q[1]$ be the 
$2$-term rigid object in $\K$ corresponding to the support $\tau$-rigid
object $C \amalg Q[1]$ in $\C(\Lambda)$.
By Lemma \ref{rigid-rigid}, we have that $\mathbb{C}_Q \amalg \PP_U$ is 
a $2$-term silting object.
\end{remark}

The next lemma, and also Lemma \ref{Bong}, appear in a more general context in \cite{bpp}.
We provide proofs adapted to our setting, for the convenience of the reader.
\begin{lemma}\label{exchange-triangle}
Let $U$ be a $\tau$-rigid module, and consider the support $\tau$-tilting object $C \amalg U \amalg Q[1]$ in $\C(\Lambda)$, where $C \amalg Q[1]$ is the co-Bongartz complement of $U$. Let $\mathbb{C}_Q = 
\PP_{C} \amalg Q[1]$ be the corresponding $2$-term rigid object in $\K$.
Let $\alpha \colon \mathbb{P}'_U \rightarrow \mathbb{C}_Q$ be a minimal right $\add  \PP_U$-approximation of $\mathbb{C}_Q$, and complete it to a triangle:
 $$\mathbb{Y} \xrightarrow{\beta} \PP'_U \xrightarrow{\alpha} \mathbb{C}_Q \to $$
in $\K$. Then the following hold.
\begin{itemize}
\item[(a)] The map $H^0(\alpha)$ is a minimal right $\add U$-approximation
of $H^0(\mathbb{C}_Q) = C$.
\item[(b)] The object $\mathbb{Y}$ is $2$-term. 
\item[(c)] The object $\mathbb{Y} \amalg \mathbb{\PP}_U$ is $2$-term silting.
\item[(d)] Let $B= H^0(\mathbb{Y})$. Then $B \amalg U$ is $\tau$-tilting, and $\mathbb{Y} =\PP_B$.
\item[(e)] The map $\beta \colon \mathbb{Y}\rightarrow \mathbb{P}'_U$ is a minimal left $\add  \PP_U$-approximation.
\end{itemize}
\end{lemma}

\begin{proof}
To prove (a), consider a map $\eta \colon U'' \to C$, with $U''$ in $\add U$.
Let $\widetilde{\eta} \colon \PP_{U''} \to \PP_{C}$ be a map such that
$H^0\left( \, \widetilde{\eta}\,\right) = \eta$. 
Consider the map
\begingroup
\renewcommand*{\arraystretch}{1.2}
\setlength\arraycolsep{10pt}
 $\gamma =\begin{pmatrix} \,\widetilde{\eta}\,\, \\0 \end{pmatrix} 
\colon \PP_{U''} \to \PP_{C} \amalg Q[1] = \mathbb{C}_Q$. 
\endgroup
Since $\alpha$ is a right $\add \PP_U$-approximation, the map 
$\gamma$ factors through $\alpha$. But then $\eta = H^0(\gamma)$
factors through $H^0(\alpha)$. 
Using Lemma \ref{lift}, we have that the minimality 
	of $\alpha$ implies that $H^0(\alpha)$ also is minimal. This gives part (a).

Since $C$ lies in $\Gen U$, the map $H^0(\alpha)$ is an epimorphism,
and part (b) then follows from Lemma \ref{exchange}(a).

Since $\mathbb{P}_U\amalg \mathbb{C}_Q$ is $2$-term silting by Remark \ref{rema}, 
part (c) now follows from Lemma \ref{exchange}(f).

We now prove part (d).
	We have that $\mathbb{Y}[1]$ is homotopic to the mapping cone of 
	$\alpha$. This mapping cone is 
$$
	\cdots \to 0 \to P_{U'}^{-1} \to P_{U'}^{0} \amalg P_{C}^{-1} \amalg Q \to P_{C}^0 \to 0 \to 
	\cdots 
	$$
	where $U'= H^0(\PP'_{U})$ is in $\add U$. 
Since the induced map
$H^0(\alpha) \colon U' \to C$ is an epimorphism, the map between the projective 
covers
$\omega \colon P_{U'}^0 \to P_{C}^0$ is a (split) epimorphism. 
Then $P_{U'}^0 = W \amalg W'$,
where $W = \ker \omega$ and $\omega$ induces an isomorphism $W'  \simeq P_{C}^0$.
We may replace the complex with an isomorphic one in which
the image of the map $ P_{U'}^{-1} \to P_{U'}^{0}$ is contained in $W$
and the second component of the map $W\amalg W' \amalg P_{C}^{-1} \amalg Q \to P_{C}^0$ is an isomorphism, with the other
components equal to zero.

Hence the complex 
$$0 \to W' \xrightarrow{\simeq} P_C^0 \to 0$$ splits off as a direct summand of the above mapping cone complex. Therefore, in the homotopy category, it is isomorphic to the complex 
\begin{equation}\label{cone}
\cdots \to 0 \to P_{U'}^{-1} \xrightarrow{\mu} W \amalg P_{C}^{-1} \amalg Q \to 0
\to  \cdots 
\end{equation}

The map $P^{-1}_{U'} \to P^{0}_{U'}$ is right minimal, that is: no direct summands of
$P^{-1}_{U'}$ are mapped to $0$. Since the image of 
the map $ P_{U'}^{-1} \to P_{U'}^{0}$ is contained in $W$, the complex (\ref{cone}) has 
no direct summands of the form $P \to 0$. 
In particular $\mathbb{Y}$  is a (minimal) projective presentation of  $\coker \mu = H^0(\mathbb{Y})$. 
Applying Lemma \ref{nokernel} to $\mathbb{Y} \amalg \PP_U$, we see 
that $H^0(\mathbb{Y} \amalg \PP_U)={B\amalg U}$ 
is a $\tau$-tilting 
object. 

Part (e) follows directly from Lemma \ref{exchange}(d).
\end{proof}

By Lemma~\ref{exchange-triangle}(d), we can write the triangle in the
statement of the lemma as:
\begin{equation}\label{big-triangle}
\PP_B \xrightarrow{\beta} \PP'_U \xrightarrow{\alpha} \mathbb{C}_Q \to 
\end{equation}

Let $\mathbb{C}_Q=\amalg_i \mathbb{X}_i$ be a decomposition of $\mathbb{C}_Q$
into indecomposable direct summands.
For each $i$, consider a
minimal right $\add \PP_U$-approximation $(\PP'_U)_i \xrightarrow{\alpha_i}  {\mathbb X}_i$. 
It is easy to check that 
$\amalg_i \alpha_i$ is a minimal right $\add \PP_U$-approximation
of $\mathbb{C}_Q$, and hence we may assume that $\alpha = \amalg_i \alpha_i$.
Hence, we obtain for each $i$ a triangle 
\begin{equation}\label{i-triangle}
\PP_{B_i} \xrightarrow{\beta_i}  (\PP'_U)_i \xrightarrow{\alpha_i}  {\mathbb X}_i \to 
\end{equation}
obtained by completing a minimal right $\add \PP_U$-approximation $\alpha_i$
of $\mathbb{X}_i$ to a triangle. 
We now have $\beta= \amalg_i \beta_i$, 
$B = \amalg_i B_i$ and $\PP_{B} = \amalg_i \PP_{B_i}$.

\begin{lemma}\label{triangle-bijection}
With notation as above, the map $\mathbb{X}_i \mapsto \PP_{B_i}$
is a bijection between the indecomposable direct summands of $\mathbb{C}_Q$ and the indecomposable direct summands of $\PP_{B}$.
\end{lemma}

\begin{proof}
By Lemma~\ref{exchange-triangle} and the above discussion, each map $\PP_{B_i} \xrightarrow{\beta_i}  (\PP'_U)_i$ is a minimal left 
$\add \PP_U$-approximation. 
\end{proof}

In particular we now have that the $B_i$ are indecomposable by Lemma~\ref{indec}.

\begin{lemma}\label{Bong}
The $\Lambda$-module $B$ is the Bongartz complement $B[U]$ of $U$; in particular $B$ is basic.
\end{lemma}

\begin{proof}
We first prove that $B$ is basic.
Suppose that $B_i \simeq B_j$ for some $i\not=j$.
Then also $\PP_{B_i} \simeq \PP_{B_j}$, and, by Lemma~\ref{triangle-bijection}, we have $\mathbb{X}_i \simeq \mathbb{X}_j$. But $\mathbb{C}_Q$ is basic, since by construction
both $C$ and $Q$ are basic (see Theorem~\ref{two-pairs}(b) and the remark afterwards). Hence $i=j$, and therefore $B$ is basic.

By Theorem \ref{two-pairs}(a), the Bongartz complement $B[U]$ of $U$ is characterized by the property that
${^{\perp}{\tau(B[U] \amalg U)}} = {^{\perp}{\tau U}}$. It is therefore sufficient to
prove that  ${^\perp{\tau U}} \subset {^\perp{\tau B}}$.
That is, we need to prove that $\Hom(X, \tau U) = 0$ implies that
$\Hom(X, \tau B) = 0$. By Lemma \ref{rigid-rigid}, this is equivalent to proving that
$\Hom_{\K}(\PP_U, \PP_X[1]) = 0$ implies that $\Hom_{\K}(\PP_B, \PP_X[1]) = 0$.

For this, consider part of the long exact sequence obtained by applying 
$\Hom_{\K}(\ , \PP_X[1])$ to the triangle (\ref{big-triangle}):

$$\Hom((\PP'_U, \PP_X[1]) \to \Hom(\PP_B, \PP_X[1]) \to \Hom(\mathbb{C}_Q, \PP_X[2]).$$
The last term vanishes since the complexes 
$\mathbb{C}_Q$ and $\PP_X$
are both $2$-term. \sloppy Hence 
$\Hom_{\K}(\PP_U, \PP_X[1]) = 0$ implies that $\Hom_{\K}(\PP_B, \PP_X[1]) = 0$,
as required.

This proves that ${^\perp{\tau U}} \subset {^\perp{\tau B}}$, and hence that ${^{\perp}{\tau(B \amalg U)}} = {^{\perp}{\tau U}}$, which implies that $B = B[U]$. 
\end{proof}

Recall that $\mathbb{C}_Q=  \PP_C \amalg Q[1]$, where $C = C[U]$ is the co-Bongartz
complement of $U$. We now focus on the indecomposable direct summands $Q'[1]$ of $\mathbb{C}_Q$, where $Q'$ is an indecomposable direct summand
of $Q$.

\begin{lemma}\label{fromproj}
Let $\mathbb{X}_i$ be an indecomposable direct summand of $\mathbb{C}_Q$ of the form $Q'[1]$, where $Q'$ is an indecomposable direct summand of $Q$.
Consider triangle (\ref{i-triangle}). The induced map 
$H^{-1}(\mathbb{X}_i) \to  H^0(\PP_{B_i})$, that is $Q' \to B_i$,
is a minimal left ${^\perp{\tau U}}$-approximation, and hence also a left
$\P({^\perp{\tau U}})$-approximation.
\end{lemma}

\begin{proof}
Let $M$ be in ${^{\perp}{\tau U}}$. Then, by Lemma \ref{rigid-rigid},
we have $\Hom_{\K}(\PP_U,\PP_M[1]) =0$. Apply $\Hom_{\K}(\ ,\PP_M)$ to the triangle 
(\ref{i-triangle}) and consider the exact sequence
$$\Hom(\PP_{B_i},\PP_M) \to \Hom(\mathbb{X}_i,\PP_M[1]) \to \Hom((\PP'_U)_i,\PP_M[1])$$
Since the last term vanishes, every map $\mathbb{X}_i[-1] \to \PP_M$ factors through
$\mathbb{X}_i[-1]\to \PP_{B_i}$.
This means that every map from $Q$ (regarded as a complex concentrated in degree $0$) to $\PP_M$ factors through
$\PP_{B_i}$. Hence, by Lemma \ref{lift}, 
the map $Q' \to B_i$ is a left ${^\perp{\tau U}}$-approximation. This map is non-zero, since
${^\perp{\tau U}}$ is sincere by \cite[Theorem 2.10]{air}.
It is therefore minimal, since $B_i$ is indecomposable.
Since $B_i$ is in
$\P({^\perp{\tau U}})$ by Lemma \ref{Bong}, the last statement also follows.
\end{proof}

Taking homology, the triangle (\ref{big-triangle}) 
induces an exact sequence
$$H^{-1}(\mathbb{C}_Q) \to H^0(\PP_B) \to H^0(\PP'_U) \to H^0(\mathbb{C}_Q) \to 0.$$
Note that $C = H^0(\mathbb{C}_Q)$ is the co-Bongartz complement of $U$ and
$B= H^0(\PP_B)$ is the Bongartz complement of $U$.
Let $U' = H^0(\PP'_U)$ and $H^{-1}(\mathbb{C}_Q) = \widetilde{Q}$, where $Q$ is a direct summand of $\widetilde{Q}$.
The sequence above becomes:
\begin{equation}\label{big-homology}
\widetilde{Q} \xrightarrow{\delta} B \xrightarrow{\mu} U' \xrightarrow{\gamma} C \to 0
\end{equation}

Since the triangles (\ref{i-triangle}) sum to the triangle (\ref{big-triangle}), we also have that the sequence
(\ref{big-homology}) is a direct sum of $n-\delta(U)$ exact sequences
\begin{equation}\label{i-homology}
Q_i \xrightarrow{\delta_i} B_i \xrightarrow{\mu_i} U_i' \xrightarrow{\gamma_i} C_i \to 0
\end{equation}
where for each $i$ either

\begin{itemize}
\item[Case (i)] $C_i$ is non-zero and is an indecomposable direct summand of $C$ (this happens when $\mathbb{X}_i$ is the minimal projective presentation of a summand $C_i$ of $C$), or
\item[Case (ii)]  $Q_i$ is an indecomposable projective direct summand of $Q$, the map $\mu_i$ is an epimorphism and $C_i=0$  (this happens when $\mathbb{X}_i$ is the shift of an indecomposable projective
direct summand of $Q$).
\end{itemize}
 
\begin{lemma}\label{sequence-bijection}
With notation as above, consider the exact sequence (\ref{i-homology}).
Then we have the following:
\begin{itemize}
\item[(a)] The map $\mu_i$ is a minimal left $\add U$-approximation.
\item[(b)] In case (i), the map $\gamma_i$ is a minimal right $\add U$-approximation, while in case (ii), the map $\gamma_i$ is the zero map and $\delta_i$ is a minimal left $\P({^\perp{\tau U}})$-approximation.
\end{itemize}
\end{lemma}

\begin{proof}
Using the fact that each $\beta_i$ in (\ref{i-triangle}) is a
minimal left $\add \PP_U$-approximation,
in combination with Lemma \ref{lift}(b), it follows that each 
$\mu_i$ is a left $\add U$-approximation. Minimality follows from the fact that 
 $\add U \cap \add C = 0$. This proves (a).

For (b), note that in case (i) the fact that $\gamma_i$ is a minimal right 
$\add U$-approximation follows from Lemma \ref{exchange-triangle}.
In case (ii), the map $\gamma_i$ must be zero as $C_i=0$.
The fact that $\delta_i$ is a minimal left $\P({^\perp{\tau U}})$-approximation
follows from Lemma~\ref{fromproj}.
\end{proof}

Let us summarize our findings.

\begin{proposition}\label{bijection-first}
Let $U$ be a $\tau$-rigid module. Let $B$ be the Bongartz complement of $U$
and $C$ the co-Bongartz complement of $U$, with corresponding support $\tau$-tilting
object $C\amalg U \amalg Q[1]$ such that $\add (C \amalg U)  =\P(\Gen U)$, as in Theorem~\ref{two-pairs}.

Then there is a triangle
\begin{equation}
\label{BCtriangle}
\PP_B \overset{\beta}{\rightarrow} \PP'_U \overset{\alpha}{\rightarrow} \mathbb{C}_Q
\rightarrow,
\end{equation}
where $\mathbb{C}_Q=\PP_C \amalg Q[1] $ and $\beta$ (respectively, $\alpha$)
is a minimal left (respectively, right) $\add \PP_U$-approximation. This triangle
is the direct sum of $n-\delta(U)$ triangles
\begin{equation}
\label{BCtrianglei}
\PP_{B_i} \overset{\beta_i}{\rightarrow} (\PP'_U)_i \overset{\alpha_i}{\rightarrow} \mathbb{X}_i \rightarrow,
\end{equation}
where $B=\amalg_i B_i$ is a decomposition of $B$ into a direct sum of indecomposable modules and the $\mathbb{X}_i$ are the indecomposable direct
summands of $\mathbb{C}_Q$.

Let $B=B'\amalg B''$, where $B'$ is the direct sum of all indecomposable
direct summands of $B$ with the property that the minimal left $\add U$-approximation $B_i\rightarrow U_i$ is not an epimorphism, and $B''$ is the complement of $B'$ in $B$.
\begin{itemize}
\item[(a)]
For each indecomposable direct summand $B_i$ of $B$ which is a summand
of $B'$, there is an exact sequence
$$B_i\xrightarrow{\mu_i} U'_i \xrightarrow{\gamma_i} C_i\to 0,$$
where $\mu_i$ (respectively, $\gamma_i$) is a \sloppy minimal left (respectively, right)
$\add U$-approximation and $C_i$ is an indecomposable direct summand of $C$.
This arises from part of the long exact sequence associated to~\eqref{BCtrianglei}:
$$H^0(\PP_{B_i})\rightarrow H^0((\PP'_U)_i)\rightarrow H^0(\mathbb{X}_i)\rightarrow H^1(\PP_{B_i})=0,$$
where $\mathbb{X}_i= \mathbb{P}_{C_i}$.
The map $B_i\mapsto C_i$ is a correspondence between the indecomposable direct summands of $B'$ and the indecomposable direct summands of $C$.

\item[(b)] For each indecomposable direct summand $B_i$ of $B$ which is a
summand of $B''$, there is an exact sequence
$$Q_i\xrightarrow{\delta_i} B_i\xrightarrow{\mu_i} U'_i\rightarrow 0,$$
where $Q_i$ is an indecomposable direct summand of $Q$, with $U'_i\in \add U$
and with $\delta_i \colon Q_i\to B_i$ a minimal left $\P({^{\perp}\tau U})=\add B$-approximation. This arises from part of the long exact sequence associated
to~\eqref{BCtrianglei}:
$$H^{-1}(\mathbb{X}_i)\rightarrow H^0(\PP_{B_i})\rightarrow H^0((\PP'_U)_i)\rightarrow H^0(\mathbb{X}_i)=0,$$
where $\mathbb{X}_i=Q_i[1]$. 
The map $B_i \mapsto Q_i$ is a bijection between the indecomposable direct summands of $B''$ and the indecomposable direct summands of $Q$.
\item[(c)] There is a bijection between the indecomposable direct summands
of $B$ and the indecomposable direct summands of the support $\tau$-rigid object $C \amalg Q[1]$ in $\C(\Lambda)$.
\end{itemize} 
\end{proposition}

\begin{proof}
Parts (a) and (b) follow from Lemmas \ref{triangle-bijection}, \ref{Bong}, \ref{sequence-bijection} and the above discussion. Part (c) is a direct consequence of parts (a) and (b).
\end{proof}

We recall the following version of Wakamatsu's lemma from \cite{air}.
\begin{lemma}\cite[Lemma 2.6]{air}\label{wakamatsu}
Let $U$ be a $\tau$-rigid module and let $\alpha \colon U' \to X$ be a right $\add U$-approximation. Then $\ker \alpha$ lies in ${^\perp(\tau U)}$.
\end{lemma}

Later we will need the following stronger version of Lemma \ref{sequence-bijection}(a), which is due to \cite{air}.

\begin{lemma}\cite[Lemma 2.20]{air}\label{Gen-app}
The map $\mu_i$ of Lemma \ref{sequence-bijection} is a minimal left $\Gen U$-approximation.
\end{lemma}

\begin{proof}
The proof is essentially contained in the proof of Lemma 2.20 of \cite{air},
but we give the details for convenience. For an object $V$ in $\Gen U$
there is a short exact sequence
\begin{equation}\label{app-seq}
0 \to Z \to U' \to V \to 0,
\end{equation}
where $U' \to V$ is a right 
$\add U$-approximation, and $Z$ lies in ${^{\perp}{\tau U}}$
by Lemma \ref{wakamatsu}. Since ${^{\perp}{\tau U}}
\subseteq {^{\perp}{\tau B}}$, we have $\Hom(Z,\tau B)= 0$, and
hence $\Ext^1(B,Z)= 0$, by the Auslander-Reiten formula.
Applying $\Hom(B_i, \ )$ to the exact sequence (\ref{app-seq})
we get an exact sequence
$$\Hom(B_i, U')  \to \Hom(B_i, V)  \to \Ext^1(B_i, Z).$$
Since the last term vanishes, the first map is an epimorphism.

Consider an arbitrary map $B_i \to V$. By the above, it factors
$B_i \to U' \to V$. But the map $\mu_i \colon B_i \to U_i$ is a
left $\add U$-approximation by Lemma~\ref{sequence-bijection}(a),
and so $B_i \to U'$ factors
$B_i \to U_i \to U'$, and hence $B_i \to V$ also factors through
$\mu_i$. This concludes the proof.
\end{proof}

\section{Reduction}\label{sec-red}
We fix a $\tau$-rigid $\Lambda$-module $U$ throughout this section.
Recall that a pair $(\T,\F)$ of subcategories of $\module \Lambda$ is called a \emph{torsion pair} if $\T = {^{\perp}\F}$ and $\F= {\T^{\perp}}$. For a given torsion pair $(\T,\F)$ and an arbitrary module $X$, there
is a (unique up to isomorphism) exact sequence $$0 \to t(X) \to X \to f(X) \to 0$$
with the property that $t(X)$ is in $\T$ and $f(X)$ is in $\F$. It is known as the \emph{canonical sequence} for $X$.

\begin{lemma}\cite[Theorem 5.10]{auslander-smalo}
The pair $(\Gen U, U^{\perp})$ is a torsion pair in $\module \Lambda$.
\end{lemma}

From now on we only consider this torsion pair, and use the notation $t$ and $f$ relative to this pair.

We next recall some results and notions, mostly from \cite{jasso}.
Let $J(U)= {U^{\perp}} \cap {^{\perp}{\tau U}}$.
The following summarizes some important facts about $J(U)$.
\begin{proposition}\label{jass}
Let $B$ be the Bongartz complement of $U$.
Then we have:
\begin{itemize}
\item[(a)]\cite[Theorem 3.8]{jasso} The subcategory $J(U)$ is equivalent to $\module \End(B \amalg U)/I$, where $I$ is the ideal generated by all maps factoring through $U$. 
\item[(b)]\cite[Cor. 3.22]{bst}, \cite[Theorem 3.28]{dirrt} The subcategory $J(U)$ is an exact abelian (wide) subcategory of $\module \Lambda$.
\item[(c)]\cite[Theorem 3.8]{jasso} If $U$ is indecomposable, then $J(U)$ has $n-1$ simple modules up to isomorphism.
\end{itemize}
\end{proposition}

The following result of Auslander and Smal{\o} is very useful.

\begin{lemma}\cite[Proposition 5.8]{auslander-smalo}\label{aus-sma}
For $X,Y$ in $\module \Lambda$ we have $\Hom(X, \tau Y) = 0$ \sloppy if and only if $\Ext^1(Y, \Gen X)= 0$
\end{lemma}

\begin{lemma}\label{Gen-lemma}
\begin{itemize}
\item[(a)] If $X$ is in $\Gen U$, then $X$ is in $\P(\Gen U)$ if and
only if $X \amalg U$ is $\tau$-rigid.
\item[(b)] If $X \amalg U$ is $\tau$-rigid, then $\Ext^1(X,t(Z))= 0$ 
for any module $Z$.
\end{itemize}
\end{lemma}

\begin{proof}
We have (see Theorem \ref{two-pairs}) that 
$\P(\Gen U) = \add (C \amalg U)$, where $C \amalg U$ is $\tau$-rigid.
Hence, if $X$ is in $\P(\Gen U)$, then $X \amalg U$ is $\tau$-rigid.

Conversely, assume $X \amalg U$ is $\tau$-rigid. Then it follows from Lemma \ref{aus-sma} that $X$ lies in $\P(\Gen U)$.
This proves (a), and (b) is a direct consequence of (a). 
\end{proof}

Our next step towards the main result, is the following bijection.

\begin{proposition}\label{object-bijection}
For a $\tau$-rigid module $U$, and with notation as before,
the map $f$ induces a bijection:  
\begin{gather*}\text{Objects $X$ in $\ind \module \Lambda$ such that $X \amalg U$ is $\tau$-rigid and $X$ is not in $\Gen U$} \\ \updownarrow \\ \text{Objects in $\ind J(U)$ which are $\tau$-rigid in $J(U)$}\end{gather*}
\end{proposition}

In order to prove Proposition \ref{object-bijection}, we will need the following lemmas.

\begin{lemma}\label{pres-ind}
If $X$ is an indecomposable $\Lambda$-module such that $X \amalg U$ is $\tau$-rigid,  then either $f(X)$ is indecomposable or $f(X)=0$.
We have $f(X) = 0$ if and only if $X$ is in $\Gen U$.  
\end{lemma}

\begin{proof}
Note that $f(X) = 0$ if and only if $X$ is in $\Gen U$, by the definition of $f$.
For the rest of the statement, it is sufficient to prove that there is a surjective ring map $\End(X) \to \End(f(X))$.
The (well-known) functoriality of $f$ gives a map from $\End(X)$ to $\End(f(X))$; we recall the construction now.
Let $\phi$ be in $\End(X)$ and
consider the diagram
$$
\xymatrix{
0 \ar[r] & t(X) \ar[r]^{\eta_X} & X \ar[r]^{u_X}  \ar[d]^{\phi} & f(X) \ar[r] & 0  \\
0 \ar[r] & t(X) \ar[r]^{\eta_X} & X \ar[r]^{u_X} & f(X) \ar[r]  & 0 
}
$$
Since $\Hom(t(X),f(X))=0$, there is a map $\theta \colon f(X) \to f(X) $ such that 
$u_X \phi  = \theta u_X$.
Since $u_X$ is an epimorphism, there is a unique map $\theta$ with this property, so this gives a well defined map $\End(X) \to \End(f(X))$. It is then easy to check
that this map is a ring map. We claim that it is surjective. 

Consider part of the long exact sequence
$$\Hom(X,X) \to \Hom(X, f(X)) \to \Ext^1(X,t(X))$$
obtained by applying $\Hom(X,\ )$ to the canonical sequence of $X$.
\sloppy We have that $\Ext^1(X,t(X))=0$ by Lemma \ref{Gen-lemma}. 
Hence the map $\Hom(X,X) \to \Hom(X, f(X))$ is surjective.
Furthermore, applying $\Hom(\ , f(X))$ to the canonical sequence gives that 
$\Hom(f(X),f(X)) \simeq \Hom(X,f(X))$, since $\Hom(t(X),f(X))= 0$. The claim follows.
\end{proof}

\begin{lemma}\label{inclusive}
Let $X,Y$ be indecomposable modules not in $\Gen U$, and such that both $U \amalg X$  and $U \amalg Y$ are $\tau$-rigid.
Then $f(X) \simeq f(Y)$ implies that $X \simeq Y$.
\end{lemma}

\begin{proof}
Fix an isomorphism $\phi \colon f(X) \to f(Y)$, and consider the diagram 
$$
\xymatrix{
0 \ar[r] & t(X) \ar[r] & X \ar[r]^{u_X} & f(X) \ar[r] \ar[d]^{\phi} & 0 \\
0 \ar[r] & t(Y) \ar[r] & Y \ar[r]^{u_Y} & f(Y) \ar[r]  & 0
}
$$
where the rows are the canonical sequences for $X$ and $Y$.
Consider part of the long exact sequence:
$$\Hom(X,Y) \xrightarrow{(X,u_Y)} \Hom(X,f(Y)) \to \Ext^1(X,t(Y)).$$
We have that $\Ext^1(X,t(Y)) =0$ by Lemma \ref{Gen-lemma}, and hence 
the map
$$\Hom(X,Y) \xrightarrow{\Hom(X,u_Y)}  \Hom(X,f(Y))$$
is surjective. Now choose a map $\psi \colon X \to Y$ satisfying $u_Y \psi  = \phi u_X$.
By a symmetric argument, we can also choose a map $\psi' \colon Y \to X$ such that $u_X \psi' =  \phi^{-1} u_Y$.

Now consider the composition $\psi' \psi$ in $\End(X)$. 
By Fitting's Lemma, any endomorphism of an indecomposable finite length module is either invertible or nilpotent.
For any positive integer $n$ we have $u_X (\psi' \psi)^n = (\phi^{-1} \phi)^n u_X = u_X \neq 0$. 
Note that $u_X$ is non-zero, since $X$ is not in $\Gen U$.
Hence $\psi' \psi$ is not nilpotent and thus an automorphism. Therefore $\psi$ is an isomorphism, and this finishes the proof.
\end{proof}

Recall from Proposition \ref{jass} that $J(U)$ is equivalent to the module category of a finite dimensional algebra which we call $\Gamma_U$.
We say that an object in $J(U)$ is \emph{$\tau$-rigid in $J(U)$} if the corresponding $\Gamma_U$-module under this equivalence is $\tau$-rigid. The Auslander-Reiten
translate in $\Gamma_U$-mod induces a translate on $J(U)$, but in general this is not the same as the Auslander-Reiten translate in $\Lambda$-mod.
In particular, a module which is $\tau$-rigid in $J(U)$ in the above sense is
not necessarily a $\tau$-rigid $\Lambda$-module (see Example~\ref{ex3}).

\begin{lemma}\label{pres-rig}
Assume $U \amalg X$ is $\tau$-rigid in $\module \Lambda$. Then either 
\begin{itemize}
\item[(i)] $X$ is in $\P(\Gen U)$ and $f(X)=0$, or 
\item[(ii)] $f(X)$ is $\tau$-rigid in $J(U)$.
\end{itemize}
\end{lemma}

\begin{proof}
Clearly $f(X) = 0$ if and only if $X$ is in $\Gen U$. 
It is also clear that $f(X)$ is in $J(U)$, since  $^{\perp}\tau U$ contains $X$ and is closed under factors. 
If $X$ is in $\Gen U$ and  $U \amalg X$ is $\tau$-rigid, then $X$ is in $\P(\Gen U)$ by Lemma \ref{Gen-lemma}.

So assume $f(X) \neq 0$. 
Consider the following exact sequence, obtained by applying $\Hom(\ ,\Gen X \cap J(U))$ to
the canonical sequence for $X$:
$$\Hom(t (X),\Gen X \cap J(U)) \to \Ext^1(f(X),\Gen X \cap J(U)) \to \Ext^1(X ,\Gen X \cap J(U))$$
Since $\Hom(X, \tau X) =0$ we have $\Ext^1(X, \Gen X) = 0$ by Lemma \ref{aus-sma}, so in particular 
$\Ext^1(X ,\Gen X \cap J(U)) = 0$. We have $\Hom(t (X),\Gen X \cap J(U)) = 0$ 
since $t(X)$ is in $\Gen U$,
and $\Hom(\Gen U, J(U)) = 0$, since $J(U) \subseteq U^{\perp}$.
Therefore also $\Ext^1(f(X),\Gen X \cap J(U))=0$.

We have $\Gen f(X) \subseteq \Gen X$, and
hence $\Gen_{J(U)} f(X) \subseteq  \Gen X \cap J(U)$, so we have
$\Ext^1(f(X),\Gen_{J(U)} f(X))=0$.
Since $J(U)$ is an extension-closed subcategory of $\module \Lambda$, it
follows from Lemma~\ref{aus-sma} that $f(X)$ is $\tau$-rigid in $J(U)$.
\end{proof}

Note that this could also be seen by completing $U \amalg X$ to a $\tau$-tilting module (using \cite{air}, and then applying~\cite[Theorem 3.15]{jasso}).

We can now complete the proof of the Proposition.

\begin{proof}[Proof of Proposition \ref{object-bijection}]
It follows from Lemma \ref{pres-rig} that $f$ sends objects $X$ in $\ind \module \Lambda$ such that 
$X \amalg U$ is $\tau$-rigid and not in $\Gen U$ to objects in $\ind J(U)$ which are $\tau$-rigid in $J(U)$.
By Lemma \ref{pres-ind} $f$ sends indecomposable modules 
to indecomposable modules and it follows from
Lemma \ref{inclusive} that $f$ induces an injective map.

It is a direct consequence of Theorem 3.15 in \cite{jasso} that 
$f$ induces a surjective map.
\end{proof}

\begin{lemma}\label{red}
The map $X \mapsto f(X)$ induces a bijection between 
$\ind \P({^{\perp} \tau U}) \setminus \ind U$ and $\ind \P(J(U))$.
\end{lemma}

\begin{proof}
By Lemma~\ref{pres-ind}, indecomposables are preserved by $f$.
The result follows from this fact and a special case of~\cite[Prop. 3.14]{jasso}.
\end{proof}

\begin{definition} \label{defpsi}
Define a map $\rho \colon
\{\ind \P(\Gen U)\setminus \ind(U) \} \cup 
\{ \ind (\P(\Lambda) \cap {^{\perp}U}) \} \to \ind \P(J(U))$
as follows.

\begin{itemize}
	\item[(i)] If $X$ lies in $\{\ind \P(\Gen U)\setminus \ind(U)  \}$,
let $\PP_{U_X} \to \PP_X$ be a minimal right $\add \PP_U$-approximation of $\PP_X$ in $D^b(\module \Lambda)$ (where $U_X$ lies in $\add U$), and complete this to a triangle $\mathbb{R}_X \to \PP_{U_X} \to \PP_X \to.$
We set $\rho(X) = f(H^0(\mathbb{R}_X))$.
	
\item [(ii)] If $X$ lies in $\ind(\P(\Lambda) \cap {^{\perp}U})$, let 
 $a_X \colon X \to B_X$ be a
	 minimal left $\P({}^\perp \tau U)$-approximation of $X$ in $\module \Lambda$,
so $B_X$ is an indecomposable direct summand of the Bongartz complement of $U$. We set $\rho(X) = f(B_X)$.

\end{itemize}
\end{definition}

\begin{remark} \label{expandrho}
Note that in case (i) we have that
since $X$ is a direct summand of the co-Bongartz complement
of $U$, it follows from Lemmas~\ref{exchange-triangle} and~\ref{Bong} that $\mathbb{R}_X=\PP_{B_X}$
for an indecomposable direct summand $B_X$ of the Bongartz complement
of $U$. By Proposition~\ref{bijection-first}(a), taking $C_i=X$, there is an exact
sequence $$B_X \xrightarrow{a_X} U_X \xrightarrow{b_X} X\rightarrow 0,$$
where $B_X=H^0(\mathbb{R}_X)$, the map $a_X$ is a minimal left $\add U$-approximation
and $b_X$ is a minimal right $\add U$-approximation.

Note that in case (ii) we have that 
since $X$ is a direct summand of $Q$, where $Q[1]$ is as in Theorem \ref{two-pairs}, it follows
from Proposition~\ref{bijection-first}(b), taking $Q_i=X$, that there is an exact sequence
$$X\xrightarrow{c_X} B_X\xrightarrow{d_X} U_X\rightarrow 0,$$
where $d_X$ is a minimal left $\add U$-approximation.
\end{remark}

Now, combining Proposition \ref{bijection-first} and Lemma \ref{red},
we obtain the following. 

\begin{proposition}\label{bijection-second}
The map $\rho$ defines a bijection 
\begin{gather*}
\{ \ind \P(\Gen U)\setminus \ind(U)  \} \cup 
\{ \ind \P(\Lambda) \cap {^{\perp}U} \} \\
\updownarrow \\
\ind \P(J(U))
\end{gather*}
\end{proposition}

We end this section with a lemma which will be useful for the proof of the main result.

Recall that a module $Y_0$ in a full subcategory $\Y$ of $\module \Lambda$ is called \emph{split projective} in $\Y$, if each epimorphism $Y \to Y_0$ with $Y$ in $\Y$ is split.

\begin{lemma}\label{split}
Let $B$ be the Bongartz complement of $U$.
Then each direct summand $B_i$ in $B$ is split projective in ${^{\perp}{\tau U}}$. 
\end{lemma}

\begin{proof}
Note that since $B$ is the Bongartz complement of $U$ we have that ${^{\perp}{\tau (B \amalg U)}} = {^{\perp}{\tau U}}$ and $\add (B \amalg U)) =\P({^{\perp} {\tau U}})$, by Theorem \ref{two-pairs}.
	
Assume $B_i$ is not split projective. Then there is an epimorphism $Q' \to B_i$, with $Q'$ in $\add B/B_i \amalg U$. This map factors through the right $\add  B/B_i \amalg U$-approximation $Q \to B_i$,
	so in particular the map $Q \to B_i$ is also an epimorphism. 
Let $B_i'$ be the kernel of $Q \to B_i$. Then, by Lemma \ref{wakamatsu},
we have that $B_i'$ is in 
${^{\perp}{(\tau (B/B_i \amalg U))}} = {^{\perp}{\tau U}}$. But then the sequence $0 \to B_i' \to Q \to B_i \to 0$ splits,
since $B_i$ is Ext-projective in ${^{\perp} {\tau U}}$. This is a contradiction.
\end{proof}

\section{Main Theorem}\label{sec-main}

In this section we complete the proof of Theorem \ref{main-int}, using the reduction technique of 
Section \ref{sec-red} and the bijection of Section \ref{sec-exc}. 
 
Note that for a indecomposable projective module $P$, there is a
primitive idempotent $e$ in $\Lambda$, such that $P \simeq \Lambda e$,
and we have that $J(P) = {P^{\perp}}$ is equivalent to 
$\module (\Lambda/ \Lambda e \Lambda)$.

\begin{lemma}
If $(\U_1, \dots, \U_{t-1}, \U_t)$ in $\C(\Lambda)$ is a
signed $\tau$-exceptional sequence, then $t \leq n$.
\end{lemma}

\begin{proof}
This is clear when $n=1$. The statement follows by induction on $n$,
since, by Proposition~\ref{jass}, we have that $\delta(\Lambda') = \delta(\Lambda)-1$, when $\module \Lambda'$ is equivalent to $J(\U_t)$.
\end{proof}

For an object $\X$ in $\C(\Lambda)$, we set 
$$\abs{\X} = \begin{cases} X, & \text{if } \mathcal{X}= X \text{ is in } \module \Lambda; \\
Y, & \text{if } \mathcal{X}= Y[1] \text{ is in } (\module \Lambda)[1]. 
\end{cases}$$
Warning: apart from $\abs{\U_t}$, the modules $|\U_i|$ arising from a signed $\tau$-exceptional sequence $(\U_1, \dots, \U_t)$ are not necessarily $\tau$-rigid in $\module \Lambda$.
However, we have the following.

\begin{lemma} \label{l:exceptional}
Let $U$ be a $\tau$-rigid object in $\module \Lambda$ and suppose that $Y$ is
$\tau$-rigid in $J(U)$. Then $\Ext^1(Y,Y)=0$.
\end{lemma}

\begin{proof}
It follows from the Auslander-Reiten formula that $\Ext^1_{J(U)}(Y,Y)=0$.
The result then follows from the fact that $J(U)$ is an extension closed subcategory 
of $\module \Lambda$.
\end{proof}

\begin{corollary} \label{c:moduleproperty}
Let $(\U_1, \dots, \U_t)$ in $\C(\Lambda)$ be a signed $\tau$-exceptional sequence.
Then $Ext^1(|U_i|,|U_i|)=0$ for all $i$.
\end{corollary}

\begin{proof}
Firstly, $\abs{\U_t}$ is $\tau$-rigid in $\module \Lambda$, giving the result for
$i=t$ using the Auslander-Reiten formula.
The module $\abs{\U_{t-1}}$ is $\tau$-rigid in $J(\U_t)$, so the result for $i=t-1$
follows from Lemma~\ref{l:exceptional}. The result for all $i$ follows from an
inductive argument.
\end{proof}

We now restate our main result.

\begin{theorem}\label{main}
Let $\Lambda$ be  a finite dimensional algebra. For each $t \in \{1, \dots, n \}$
there is a bijection between the set of ordered support $\tau$-rigid
objects of length $t$ in $\C(\Lambda)$ and the
set of signed $\tau$-exceptional sequences of length $t$ in $\C(\Lambda)$. 
\end{theorem}

For $t=n$, we obtain.

\begin{corollary}\label{main-cor}
Let $\Lambda$ be  a finite dimensional algebra. 
There is a bijection between the set of ordered support $\tau$-tilting
objects in $\C(\Lambda)$ and the
set of complete signed $\tau$-exceptional sequences in $\C(\Lambda)$. 
\end{corollary}

The remainder of this section is devoted to proving Theorem \ref{main}.
The main idea of the proof is to work by induction on $t$, making use of
Propositions \ref{object-bijection} and \ref{bijection-second}, which we for convenience reformulate below.

Recall the definition of $\rho$ from Definition~\ref{defpsi} (see also Proposition~\ref{bijection-second}). Also, note that by Lemma \ref{Gen-lemma}, we have that for $U$ $\tau$-rigid and $X \in \Gen U$, then $X \amalg U$ is $\tau$-rigid if and only if $X$ is in $\P(\Gen U)$. Hence, if $X\in\ind \Gen U\setminus \ind U$ and $X\amalg U$ is $\tau$-rigid, then $\rho$ can be applied to $X$.

\begin{proposition}\label{third-bijection}
Let $U$ be a $\tau$-rigid module in $\module \Lambda$.
Then there is a bijection $\E_U$ between the sets
$$\{ X \in \ind \module \Lambda \setminus \ind U \mid
 X \amalg U \text{ is $\tau$-rigid} \} \cup 
\{ \ind(\P(\Lambda) \cap {^{\perp}U})[1] \}$$
and
$$\{X \in \ind J(U) \mid 
 X \text{ $\tau$-rigid in $J(U)$} \} \cup 
\{ \ind(\P(J(U))[1] \}$$
given by
$$ \E_{U}(X) = 
\begin{cases}
f(X), & X \in \{ \ind \module \Lambda \setminus \ind  U \mid
 X \amalg U \text{ is $\tau$-rigid and } X \not \in \Gen U \}; \\
\rho(X)[1], & \parbox[t]{.75\textwidth}{$X \in \{ \ind \module \Lambda \setminus \ind U \mid X \amalg U$ is $\tau$-rigid and $X \in \Gen U \};$}\\
\rho(X[-1])[1],  & X\in \ind(\P(\Lambda) \cap {^{\perp}U})[1].
\end{cases} 
$$

\end{proposition}

\begin{remark}
The reduction technique of Proposition~\ref{third-bijection} can be seen as an extension of Jasso's reduction technique \cite[Cor.\ 3.18]{jasso}. Indeed, Jasso's technique forms an important part of the proof.
However, it is important to note that modules generated by $U$ (using the notation of Proposition~\ref{third-bijection}) (which are killed by Jasso's map) are not sent to modules, but rather to
`signed' objects $P[1]$ in $\C(J(U))$ (where $P$ is a relative projective in $J(U)$). The cases of our proof dealing with such objects do not follow from Jasso's theorem.

The bijection of Proposition~\ref{third-bijection} is needed for the proof of our main result,
Theorem~\ref{main}, even in the complete case.

We note that in recent work~\cite{borve} (after this paper appeared on the Arxiv) B\o rve used the correspondence between support $\tau$-tilting modules and $2$-term silting objects, and results from \cite[Section 4]{jasso} to give an alternative proof of Proposition~\ref{third-bijection} using
Iyama-Yoshino reduction \cite[4.2, 4.7]{iy}.
\end{remark}

We extend the domain of $\E_U$ to $$\add (\{ X \in \ind \module \Lambda \setminus \ind U \mid
 X \amalg U \text{ is $\tau$-rigid} \} \cup 
\{ \ind(\P(\Lambda) \cap {^{\perp}U})[1] \})$$
by setting $\E_U(\amalg_{i=1}^s X_i) = \amalg_{i=1}^s \E_U(X_i)$
for objects $X_1,X_2,\ldots X_s$ in
$$\{ X \in \ind \module \Lambda  \setminus \ind U \mid
 X \amalg U \text{ is $\tau$-rigid} \} \cup 
\{ \ind(\P(\Lambda) \cap {^{\perp}U})[1] \}.$$

\begin{proposition} \label{strongthirdbijection}
Let $U$ be a $\tau$-rigid module in $\module \Lambda$ with $\delta(U) = t'$.
For any positive integer $t \leq n- t'$, 
the map $\E_U$ induces a bijection between the set of support $\tau$-rigid objects $X$ in $\C(\Lambda)$ such that
$\delta(X) = t$, $X \amalg U$ is support $\tau$-rigid and $\add X \cap \add U=0$, and the set of support $\tau$-rigid objects $Y$ in $\C(J(U))$ with $\delta(Y) = t$.
\end{proposition}

\begin{proof}
We first need to prove that for any $X$ in $\C(\Lambda)$, if $X \amalg U$ is support $\tau$-rigid
with $\add X\cap \add U=0$, then
$\E_U(X) = \widetilde{X}$ is support $\tau$-rigid in $\C(J(U))$.
For this, let $X_i$ and $X_j$ be two indecomposable direct summands in $X$, and consider the
following cases.

\bigskip

\noindent {\bf Case I:}
Let $X_i,X_j$ both be in $\module \Lambda$ and not in $\Gen U$.
Then by Lemma \ref{pres-rig}, we have that 
$$\widetilde{X}_i \amalg \widetilde{X}_j =f(X_i) \amalg f(X_j) =f(X_i \amalg X_j)$$
is $\tau$-rigid in $J(U)$. 

\bigskip

\noindent {\bf Case II:}
Let $X_i$ be in $\Gen U$, and assume $X_j$ is in $\module \Lambda$ but not in $\Gen U$.

Then $\widetilde{X_i}[-1]$ is in $\P(J(U))$, and we need to prove 
that $\Hom(\widetilde{X}_i[-1],\widetilde{X}_j) = 0$.
Note that
$$\widetilde{X}_i=\E_U(X_i)=\rho(X_i)[1]=f(H^0(\mathbb{R}_{X_i}))[1]=
f(B_{X_i})[1]$$
(see Definition \ref{defpsi}, Remark~\ref{expandrho} and Proposition \ref{third-bijection}),
where there is an exact sequence 
$$B_{X_i}\xrightarrow{a_{X_i}} U_{X_i} \rightarrow X_i\rightarrow 0$$
with $a_{X_i}$ a minimal left $\add U$-approximation.

Moreover $\widetilde{X}_j = f(X_j)$, so we need to prove that
$\Hom(f(B_{X_i}),f(X_j)) = 0$.
We have that $B_{X_i}$ is in $\P({^{\perp}{\tau U}})$ by Lemma \ref{Bong}, and
$f(B_{X_i})$ is in $\P(J(U))$ by Lemma \ref{red}.

We have that $0 \to  \Hom(f(B_{X_i}),f(X_j)) \to \Hom(B_{X_i},f(X_j))$ is exact,
so it suffices to show that $\Hom(B_{X_i},f(X_j))= 0$. 
For this, apply $\Hom(B_{X_i}, \ )$ to the exact sequence
$$0 \to t(X_j) \to X_j \to f(X_j) \to 0$$
and consider the long exact sequence
$$0 \to  \Hom(B_{X_i},t(X_j)) \to  \Hom(B_{X_i},X_j) \to  \Hom(B_{X_i},f(X_j)) \to \Ext^1(B_{X_i},t(X_j)).$$
We have that $t(X_j)$ is in $\Gen U \subseteq  {^{\perp}{\tau U}}$, and
$B_{X_i}$ is in $\P( {^{\perp}{\tau U}})$, so $\Ext^1(B_{X_i},t(X_j)) =0$.
Therefore it is sufficient to show that $\Hom(B_{X_i},t(X_j)) \to  \Hom(B_{X_i},X_j)$
is an epimorphism.

So consider an arbitrary map $g \colon B_{X_i} \to X_j$. We first
claim that this map factors through an object in $\add U$. 
Recall from Definition~\ref{defpsi} that there is a triangle
$$\PP_{B_{X_i}} \to \PP_{U_{X_i}} \to \PP_{X_i} \to $$
with $U_{X_i}$ in $\add U$.
Now apply $\Hom_{\K}(\ , \PP_{X_j})$ to this triangle, and consider the
\sloppy exact sequence
$$\Hom(\PP_{U_{X_i}},\PP_{X_j}) \to \Hom(\PP_{B_{X_i}},\PP_{X_j}) \to \Hom(\PP_{X_i},\PP_{X_j}[1])$$
By assumption we have $\Hom(X_j, \tau X_i) =0$, and this implies
that the last term $\Hom(\PP_{X_i},\PP_{X_j}[1])$ vanishes, 
by Lemma \ref{rigid-rigid}.
Hence the map $\Hom(\PP_{U_{X_i}},\PP_{X_j}) \to \Hom(\PP_{B_{X_i}},\PP_{X_j})$ is an epimorphism, and so any map 
$\PP_{B_{X_i}} \to \PP_{X_j}$ factors through $\PP_{U_{X_i}}$. By Lemma \ref{lift},
this means that the map $g$ factors through $U_{X_i}$ in $\add U$.
Assume $B_{X_i} \xrightarrow{g} X_j = B_{X_i} \xrightarrow{h} U_{X_i} \xrightarrow{p} X_j$.
Then $U_{X_i} \xrightarrow{p} X_j$ factors through $t(X_j) \to X_j$. Hence $g$ also factors through  $t(X_j) \to X_j$, and
we have proved the claim that  $\Hom(B_{X_i},t(X_j)) \to  \Hom(B_{X_i},X_j)$
is an epimorphism. It then follows that 
$\Hom(\widetilde{X}_i[-1],\widetilde{X}_j) = \Hom(B_{X_i}, f(X_j)) = 0$.

\bigskip

\noindent {\bf Case III:}
Now assume $X_i =P[1]$ for $P$ in $\P(\Lambda) \cap {^{\perp}U}$, and 
$X_j$ lies in $\module \Lambda$ and not in $\Gen U$. Then
$\Hom(P, X_j)= 0$.
We need to prove that $\Hom(\widetilde{X}_i[-1],\widetilde{X}_j)= 0$.
By Definition~\ref{defpsi}, there is an exact sequence:
\begin{equation}
\label{seq}
P\xrightarrow{c_{P}} B_{P}\xrightarrow{d_P} U_{P}\rightarrow 0
\end{equation}

By Definition~\ref{defpsi}, we have $\widetilde{X_i}=f(B_P)[1]$ and $\widetilde{X}_j = f (X_j)$.

Note that $\Hom(P, X_j) = 0$ implies that $\Hom(P, \widetilde{X}_j) = 0$,
since $P$ is projective. We also have $\Hom(U, \widetilde{X}_j) =
\Hom(U, f(X_j)) = 0$, by definition of $f$.

Applying $\Hom(\ , \widetilde{X}_j)$ to (\ref{seq}) gives an exact sequence
$$0\to \Hom(U_P, \widetilde{X}_j) \to \Hom(B_P, \widetilde{X}_j)  \to 
\Hom(P, \widetilde{X_j})$$
The end terms vanish, and hence we obtain that $\Hom(B_P, \widetilde{X}_j) = 0$. 
Clearly then also $\Hom(\widetilde{X}_i[-1], \widetilde{X}_j) = \Hom(f(B_P), \widetilde{X}_j) = 0$.
This finishes the proof of Case III.
Now combining the Cases I,II and III, the claim that 
$\E_U(X)$ is a support $\tau$-rigid object in $J(U)$ follows.

Now, let $M=\amalg_{i=1}^t M_i$ be a support $\tau$-rigid object in $\C(J(U))$ with $\delta(M) = t$, where each $M_i$ is indecomposable.
There are, by Proposition \ref{third-bijection}, unique indecomposable modules
$X_i$ such that $\E_U(X_i) = M_i$, the object $U \amalg X_i$ is support $\tau$-rigid,
and $X_i\not\in \add U$.
We need to prove that 
$X = \amalg_{i=1}^t X_i$ is support $\tau$-rigid as well.
For this we consider two arbitrary summands $X_i, X_j$ in $X$ and the following cases.

\bigskip

\noindent {\bf Case I:}
Assume $X_i = P[1]$ with $P$ in $\P(\Lambda) \cap {^{\perp}U}$, and assume 
that $X_j$ is in $\module \Lambda$ is such that $X_j \amalg U$ is $\tau$-rigid, $X_j$ is not in $\Gen U$,
and $\Hom(\widetilde{X}_i[-1],\widetilde{X}_j)= 0$.
We need to prove that $\Hom(P,X_j)= 0$.

Note that by Definition~\ref{defpsi},
we have $\widetilde{X}_i[-1] = f (B_P)$, where
$P \to B_P$ is a minimal $\P({^{\perp}{\tau U}})$-approximation of $P$.
Moreover $\widetilde{X}_j = f(X_j)$. 

We have that $\Hom(f(B_P), f(X_j)) = \Hom(\widetilde{X}_i[-1],\widetilde{X}_j)= 0$.
Apply
$\Hom(\ ,f(X_j))$ to the exact sequence $0 \to t(B_P) \to B_P \to f(B_P) \to 0$
and consider the exact sequence 
$$0 \to \Hom(f(B_P), f(X_j)) \to \Hom(B_P, f(X_j)) \to \Hom(t(B_P), f(X_j)) $$
The last term vanishes since $t(B_P)$ is in $\Gen U$ and $f(X_j)$ is in $U^{\perp}$.
We then obtain that also  
$\Hom(B_P,f(X_j))= 0$.

Since by assumption $X$ is in ${^{\perp}{\tau U}}$, we also have
that $f(X)$ is in ${^{\perp}{\tau U}}$. Hence any map $P \to f(X_j)$ factors
through $B_P \to f(X_j)$, since $P \to B_P$ is a left
${^{\perp}{\tau U}}$-approximation (see Lemma \ref{fromproj}). Since $\Hom(B_P,f(X_j))= 0$, we have 
$\Hom(P,f(X_j))= 0$.
We also have that $\Hom(P,U) = 0$ implies that $\Hom(P, \Gen U) = 0$, so in
particular $\Hom(P,t(X_j)) = 0$. 

Since  $\Hom(P,f(X_j)) = 0 = \Hom(P,t(X_j))$, we indeed also 
have  $\Hom(P,X_j) = 0$, and this finishes Case I.

\bigskip

\noindent {\bf Case II:} 
If $X_i, X_j$ are both in $\module \Lambda$ and not in $\Gen U$,
then $\widetilde{X}_i = f(X_i), \widetilde{X}_j =f(X_j)$ is such that 
$\widetilde{X}_i \amalg \widetilde{X}_j = \widetilde{X_i \amalg X_j}$
is $\tau$-rigid in $J(U)$, then $X_i \amalg X_j$ is $\tau$-rigid according to 
\cite[Cor. 3.18]{jasso}.

\bigskip

\noindent {\bf Case III:}
Now assume $X_i$ is in $\P(\Gen U)$, and hence $\widetilde{X}_i$ in
$\P(J(U))[1]$, while $X_j$ is in $\module \Lambda$ but not in $\Gen U$.

By assumption $\Hom(\widetilde{X}_i[-1], \widetilde{X}_j)= 0$,
and we need to prove that $\Hom(X_i, \tau X_j)= 0 = \Hom(X_j, \tau X_i)$
 
Since $X_i$ is in $\Gen U$, and $\Hom(U,\tau X_j)= 0$, we have that also
$\Hom(X_i, \tau X_j) = 0$.

By Lemma \ref{rigid-rigid} we have that in order to prove that  
$\Hom(X_j, \tau X_i) = 0$, it is sufficient to prove that  
$\Hom_{\K}(\PP_{X_i}, \PP_{X_j}[1]) = 0$.

We apply $\Hom_{{\K}}(\ , \PP_{X_j})$ to the triangle 
$$\PP_{B_{X_i}} \to \PP_{U_{X_i}} \to \PP_{X_i} \to $$
(see Definition~\ref{defpsi})
and consider the exact sequence
\begin{equation}\label{ses0}
\Hom(\PP_{U_{X_i}},\PP_{X_j}) \to \Hom(\PP_{B_{X_i}},\PP_{X_j}) \to \Hom(\PP_{X_i},\PP_{X_j}[1]) \to \Hom(\PP_{U_{X_i}}, \PP_{X_j}[1])
\end{equation}
We first note that $\Hom(X_j, \tau U)= 0$ implies that the last term
$\Hom_{\K}(\PP_{U_{X_i}}, \PP_{X_j}[1])$ vanishes.
It is therefore sufficient to prove that the first map
$\Hom(\PP_{U_{X_i}},\PP_{X_j}) \to \Hom(\PP_{B_{X_i}},\PP_{X_j})$
is an epimorphism, that is: we claim that any map $\PP_{B_{X_i}} \to \PP_{X_j}$ 
factors through $\PP_{B_{X_i}} \to \PP_{U_{X_i}}$. For this, it is sufficient that
any map $B_{X_i} \to X_j$ factors through $B_{X_i} \to U_{X_i}$. Consider the exact
sequence
\begin{equation}\label{ses1}
\Hom(B_{X_i},t(X_j)) \to \Hom(B_{X_i},X_j) \to \Hom(B_{X_i},f(X_j))
\end{equation}
obtained by applying $\Hom(B_{X_i},\ )$ to the canonical sequence for $X_j$.
We claim that the last term vanishes. 
For this consider the exact sequence 
\begin{equation}\label{ses2}
\Hom(f(B_{X_i}),f(X_j)) \to \Hom(B_{X_i},f(X_j)) \to \Hom(t(B_{X_i}),f(X_j))
\end{equation}
obtained by applying $\Hom(\ ,f(X_j))$ to the canonical sequence for $B_{X_i}$.
In (\ref{ses2}) the first term vanishes by assumption, and the last term
vanishes since $t(B_{X_i})$ is in $\Gen U$ and $f(X_j)$ is in $U^{\perp}$.
Hence $\Hom(B_{X_i},f(X_j))$, which is the last term of sequence (\ref{ses1})
also vanishes. This means that any map $B_{X_i} \to X_j$ factors 
$B_{X_i} \to t(X_j) \to X_j$. 

By Lemma \ref{Gen-app}, the map $B_{X_i} \to U_{X_i}$ is a $\Gen U$-approximation.
So the map  $B_{X_i} \to t(X_j)$ factors  $B_{X_i} \to U_{X_i} \to t(X_j)$, and hence
our original map $B_{X_i}  \to X_j$  factors $B_{X_i} \to U_{X_i} \to t(X_j) \to X_j$.
We have now proved that any map $B_{X_i} \to X_j$ factors through $B_{X_i} \to U_{X_i}$,
and hence any map $\PP_{B_{X_i}} \to \PP_{X_j}$ 
factors through $\PP_{B_{X_i}} \to \PP_{U_{X_i}}$. Therefore
$\Hom(\PP_{X_i},\PP_{X_j}[1])=0$, which implies $\Hom(X_j,\tau X_i) =0$,
and the claim is proved. This finishes the proof of case III.

Combining Cases I, II and III, and noting that the remaining cases can be checked trivially, proves the claim that $X = \amalg X_i$ is 
support $\tau$-rigid, and this concludes the proof of the proposition.
\end{proof}

We have the following special case of Proposition~\ref{strongthirdbijection}.

\begin{corollary}\label{reduction-one}
Let $t>1$. Let $U$ be an indecomposable $\tau$-rigid module.
The map $\E_U$ induces a bijection between
ordered support $\tau$-rigid objects in $\C(\Lambda)$ with last term $U$ and
length $t$, and ordered support $\tau$-rigid objects of length $t-1$ in
$\C(J(U))$. 
\end{corollary}

We need also to deal with the case where the last term in an ordered support $\tau$-rigid objects in $\C(\Lambda)$
is of the form $P[1]$. For this, we first observe the following.

\begin{lemma}\label{proj-bi}
Let $P$ be a projective $\Lambda$-module, and consider $J(P[1]) = P^{\perp}$.
\begin{itemize}
\item[(a)] The $\tau$-rigid modules in $J(P[1])$ are exactly the $\tau$-rigid modules $X$ in $\module \Lambda$
with $\Hom(P,X) =0$.
\item[(b)] The map $Q \mapsto f_P(Q)$ gives a bijection between the indecomposables in $\P(\Lambda) \setminus
\add P$ and the indecomposables in $\P(P^{\perp})$.
\end{itemize}
\end{lemma}

\begin{proof}
Part (a) is a direct consequence of \cite[Lemma 2.1]{air}.
Part (b) follows from Lemma \ref{red}. 
\end{proof}

For a projective module $P$ in $\P(\Lambda)$, consider the map
$\E_{P[1]}$ from $$\add \{ X \in \ind \module \Lambda \mid X \text{  } \tau \text{-rigid, and }  
\Hom(P,X)= 0 \} \cup
(\ind \P(\Lambda) \setminus \ind P)[1]$$ to 
$$\add \{ X \in \ind J(P[1])) \mid X \text{ }  \tau\text{-rigid} \} \cup
\ind \P(J(P[1]))[1]$$ 
defined as follows. For $X$ indecomposable, we set
$$\E_{P[1]} (X)= \begin{cases}
X & \text{if } X \text{ is } \tau\text{-rigid, and }
\Hom(P,X)= 0 \\
f_P(X)[1]  & \text{if } X \in \P(\Lambda) \setminus \ind P
\end{cases}$$
For $X = \amalg_{i=1}^t X_i$, with each $X_i$ in 
$$\{ X \in \ind \module \Lambda \mid X \text{  } \tau\text{-rigid, and } 
\Hom(P,X)= 0 \} \cup
(\ind \P(\Lambda) \setminus \ind P)[1],$$
we set $\E_{P[1]}(X) = \amalg_{i=1}^t \E_{P[1]}(X_i)$.

\begin{proposition}\label{yab}
Let $P$ be in $\P(\Lambda)$ with $\delta(P) = t'$.
\begin{itemize}
\item[(a)] 
The map $\E_{P[1]}$
is a bijection between the sets
$$\{ X \in \ind \module \Lambda  \mid
 X \text{ is $\tau$-rigid and } \Hom(P,X)=0 \} \cup 
\{ \ind(\P(\Lambda)) \setminus \ind P)[1] \}$$
and
$$\{X \in \ind J(P[1]) \mid 
 X \text{ is $\tau$-rigid} \} \cup 
\{ \ind \P(J(P[1]))[1] \}.$$
\item[(b)] For any positive integer $t \leq n-t'$, the map 
$\E_{P[1]}$ induces a bijection between the set of support $\tau$-rigid
objects $X$ in $\C(\Lambda)$ such that $\delta(X) = t$, the object
$X \amalg P[1]$ is support $\tau$-rigid and $\add X\cap \add P[1]$=0,
and the set of support  $\tau$-rigid objects in $\C(P^{\perp})$ with $t$ indecomposable direct summands.
\end{itemize}
\end{proposition}

\begin{proof}
Part (a) follows directly from Lemma \ref{proj-bi}.

\bigskip

\noindent Part (b): 
Let $Q\not\in\add P$ be an indecomposable module in $\P(\Lambda)$, and let $X$ be in $P^{\perp}$.
Apply $\Hom(\ ,X)$ to the canonical sequence
$$0 \to t_P (Q) \to Q \to f_P (Q) \to 0.$$ 
Since $\Hom(P,X)=0$, we have that $\Hom(\Gen P,X)=0$ and thus
that $\Hom(t_P (Q), X)= 0$.
It follows that $\Hom(Q, X) \simeq \Hom(f_p (Q), X)$.
The claim follows from combining this with part (a).
\end{proof}

\begin{corollary}\label{reduction-two}
Let $t>1$. Let $P$ be an indecomposable projective module.
Then the map $\E_{P[1]}$ induces a bijection between ordered
support $\tau$-rigid objects in $\C(\Lambda)$ with last term $P[1]$ and length $t$,
and ordered support $\tau$-rigid objects of length $t-1$ in $\C(J(P[1]))=\C(P^{\perp})$.
\end{corollary}

We can now prove Theorem \ref{main}.

\begin{proof}[Proof of Theorem \ref{main}]
We prove the result by induction on $\delta(\Lambda)=n$.
The statement for $t=1$ is clear. In particular, this deals with the case $n=1$.
So we assume the result to be true for algebras with a smaller number of
indecomposable projective modules up to isomorphism.
Let $U$ be an indecomposable $\tau$-rigid object in $\module \Lambda$.
By Corollary \ref{reduction-one}, there is a bijection between
the ordered support $\tau$-rigid objects in $\C(\Lambda)$ ending  in $U$ and the ordered support $\tau$-rigid objects in $\C(J(U))$.
The ordered support $\tau$-rigid objects in $\C(J(U))$ are by the induction hypothesis in bijection with the signed $\tau$-exceptional sequences in $\C(J(U))$.
And by definition a sequence $(U_1, \dots, U_{t-1})$ is a signed $\tau$-exceptional sequence in $\C(J(U))$ if and only if 
$(\U_1, \dots, \U_{t-1}, U)$ is a signed $\tau$-exceptional sequence in $\C(\Lambda)$.

Let $P$ be an indecomposable module in $\P(\Lambda)$. In a similar way to the above, there is a bijection between ordered support $\tau$-rigid objects in $\C(\Lambda)$
ending in $P[1]$ and signed $\tau$-exceptional sequences in $\C(J(P[1]))$,
using the induction hypothesis and Corollary \ref{reduction-two}.
\end{proof}

\begin{remark}
We now give a more explicit description of the bijection constructed in 
Theorem \ref{main}, and the inverse of this bijection.
Let $\W$ denote a wide subcategory of $\module \Lambda$ which is equivalent
to a module category, and let $n_{\W}$ denote its rank.
Let $\U$ denote an indecomposable $\tau$-rigid object in $\C(\W)$, that is 
either $\U = U$ for an indecomposable $\tau$-rigid module in $\W$, or $\U = P[1]$, where
$P$ is indecomposable projective in $\W$. Recall, in particular, that by 
Proposition \ref{jass},
we have that $J_{\W}(\U)$ is a wide subcategory of $\W$, and 
hence of $\module \Lambda$, equivalent to a module category.

Now, consider the bijections obtained
by combining Propositions \ref{third-bijection} and \ref{yab}.

$$\{ \mathcal{X} \in \ind \C(\W) \setminus \U \mid
 \mathcal{X} \amalg \U \text{ is $\tau$-rigid} \}$$
$$\E_{\U}^{\W}  \downarrow \text{  } \uparrow \F^{\W}_{\U}$$
$$\{\mathcal{X} \in \ind \C(J(\U)) \mid 
 \mathcal{X} \text{ $\tau$-rigid in $J(\U)$} \}$$
 
Consider, for each 
$t = 1, \dots, n_{\W}$, the bijections
$$\{\text{ordered $\tau$-rigid objects in $\W$ with $t$ indecomposable direct summands} \}$$
$$\Psi^{\W}_t  \downarrow \text{  } \uparrow \Phi^{\W}_t$$
$$\{\tau \text{-exceptional sequences in $\W$ of length $t$} \}$$
where $\Psi^{\W}_t$ is the bijection constructed in Theorem \ref{main},
and $\Phi^{\W}_t$ is its inverse.
Then we have 
$$\Psi^{\W}_t (\T_1, \dots, \T_t) = (\Psi_{t-1}^{J_{\W}(\T_t)}(\E^{\W}_{\T_t}(\T_1), \dots, \E^{\W}_{\T_t}(\T_{t-1})), \T_t)$$

Now let 
\begin{align*}
\W_t &= \W  & \U_t &= \T_t\\
\W_{t-1} &= J_{\W_t}(\U_t) &   \U_{t-1} &= \E_{\U_t}^{\W_t}(\T_{t-1})\\
& \; \; \vdots & & \; \;\vdots\\
\W_i &= J_{\W_{i+1}}(\U_{i+1}) & \U_i &=  \E_{\U_{i+1}}^{\W_{i+1}} \dots \E_{\U_{t-1}}^{\W_{t-1}} \E_{\U_t}^{\W_{t}}(\T_{i}) \\
& \; \;\vdots & & \; \;\vdots\\
\W_1 &= J_{\W_{2}}(\U_{2}) & \U_1 &=  \E_{\U_{2}}^{\W_{2}} \dots \E_{\U_{t-1}}^{\W_{t-1}} \E_{\U_t}^{\W_{t}}(\T_{1}) 
\end{align*}

It is then straightforward to verify that 
$$\Psi^{\W}_t (\T_1, \dots, \T_t) = (\U_1, \dots, \U_t).$$
and that the inverse bijection  
is given by
$$\Phi_t^{\W} (\U_1,\ldots ,\U_t)= 
(\F_{\U_t}^{\W_t}\cdots \F_{\U_2}^{\W_2}(\U_1),
\F_{\U_t}^{\W_t}\cdots \F_{\U_3}^{\W_3}(\U_2),
\cdots ,\U_t)$$
where $\W_t=\W$ and $\W_i=J_{\W_{i+1}}(\U_{i+1})$ for all $i$, as above.

\end{remark}

\section{Examples}\label{sec-ex}
Each example is given as the path algebra of a quiver modulo an admissible ideal of relations generated by paths. For each vertex $i$ of the quiver, we denote by $P_i, I_i$, $S_i$ the corresponding indecomposable projective (respectively, indecomposable injective, simple) module.

\subsection{Example 1:}
Let $Q$ be the quiver $1 \longrightarrow 2$, and let $\Lambda = k Q$.
There are three indecomposable modules,
$P_1, P_2 = S_2, S_1$, and the AR-quiver of $\module \Lambda$ is:
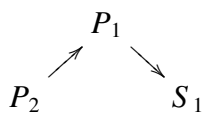
$$
\xymatrix@=4mm{
& P_1 \ar[dr] & \\
P_2 \ar[ur] \ar@{..}[rr] & & S_1
}
$$
There are 5 support $\tau$-tilting modules (= support tilting modules, since
$\Lambda$ is hereditary), and hence 10 ordered support $\tau$-tilting modules.
We list these in the table below, together with the corresponding 
complete signed $\tau$-exceptional sequences.
\begin{center}
\setlength{\extrarowheight}{1.2mm}
\begin{tabular}{|c|c||c|c|}
\hline 
\begin{tabular}{c} Ordered support \\ $\tau$-tilting object \end{tabular} &
\begin{tabular}{c} Signed $\tau$-exc. \\ sequence \end{tabular} &
\begin{tabular}{c} Ordered support \\ $\tau$-tilting object \end{tabular} &
\begin{tabular}{c} Signed $\tau$-exc. \\ sequence \end{tabular}
\\ \hline 
$(P_2, P_1)$ & $(P_2, P_1)$ & $(P_2, P_1[1])$ & $(P_2, P_1[1])$ \\ \hline
$(S_1, P_1)$ & $(P_2[1], P_1)$ & $(P_2[1], P_1[1])$ & $(P_2[1], P_1[1])$ \\ \hline
$(P_1, P_2)$ & $(S_1, P_2)$ & $(S_1, P_2[1])$ & $(S_1, P_2[1])$ \\ \hline
$(P_1[1], P_2)$ & $(S_1[1], P_2)$ & $(P_1[1], P_2[1])$ & $(S_1[1], P_2[1])$ \\ \hline  
$(P_1, S_1)$ & $(P_1, S_1)$ & &  \\ \hline
$(P_2[1], S_1)$ & $(P_1[1], S_1)$ & & \\ \hline
\end{tabular} 
\end{center}

\subsection{Example 2:}
Let $Q'$ be the quiver 
$\xymatrix{ 1 \ar^{\alpha}@<3pt>[r] & 2 \ar^{\beta}@<3pt>[l]}$,  and
let $\Lambda' = kQ'/I$, where $I$ is the ideal generated by the path 
$\beta \alpha$.
There are 5 indecomposable modules,
and the AR-quiver is:
$$
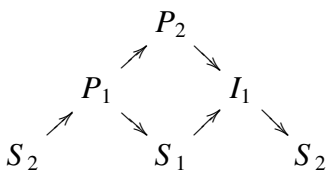
\xymatrix@=3mm{
& & P_2  \ar[dr] & & \\
& P_1 \ar@{..}[rr] \ar[dr] \ar[ur]&  & I_1  \ar[dr] & \\
S_2 \ar@{..}[rr] \ar[ur] & & S_1 \ar@{..}[rr]  \ar[ur] &  & S_2
}
$$
Note that the module $I_1$ is {\em not} $\tau$-rigid in $\module \Lambda$, while
the other four indecomposable modules are $\tau$-rigid.

There are 6 support $\tau$-tilting modules, and hence 12 ordered support $\tau$-tilting modules
We list these in the table below, together with the corresponding 
complete signed $\tau$-exceptional sequences.
\begin{center}
\setlength{\extrarowheight}{1.2mm}
\begin{tabular}{|c|c||c|c|}
\hline 
\begin{tabular}{c} Ordered support \\ $\tau$-tilting object \end{tabular} &
\begin{tabular}{c} Signed $\tau$-exc. \\ sequence \end{tabular} &
\begin{tabular}{c} Ordered support \\ $\tau$-tilting object \end{tabular} &
\begin{tabular}{c} Signed $\tau$-exc. \\ sequence \end{tabular}
\\ \hline 
$(P_1, P_2)$ & $(S_1, P_2)$ & $(P_2[1], S_1)$ & $(P_1[1], S_1)$ \\ \hline
$(S_2, P_2)$ & $(S_1[1], P_2)$ & $(P_1, S_1)$ & $(P_1, S_1)$ \\ \hline
$(P_2, P_1)$ & $(S_2, P_1)$ & $(S_2, P_1[1])$ & $(S_2, P_1[1])$ \\ \hline
$(S_1, P_1)$ & $(S_2[1], P_1)$ & $(P_2[1], P_1[1])$ & $(S_2[1], P_1[1])$ \\ \hline  
$(P_2, S_2)$ & $(I_1, S_2)$ &  $(S_1, P_2[1])$& $(S_1, P_2[1])$  \\ \hline
$(P_1[1], S_2)$ & $(I_1[1], S_2)$ &  $(P_1[1], P_2[1])$ &  $(S_1[1], P_2[1])$\\ \hline
\end{tabular} 
\end{center}

\subsection{Example 3:}\label{ex3}
Let $Q''$ be the quiver 
$$
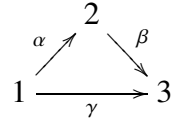
\xymatrix@=5mm{
& 2  \ar^{\beta}[dr] & \\
1  \ar^{\alpha}[ur]  \ar_{\gamma}[rr]& & 3
} 
$$
and let $\Lambda'' = kQ''/I$ where $I$ is the ideal generated by the path 
$\beta \alpha$. There are 9 indecomposable modules, and the AR-quiver is
$$
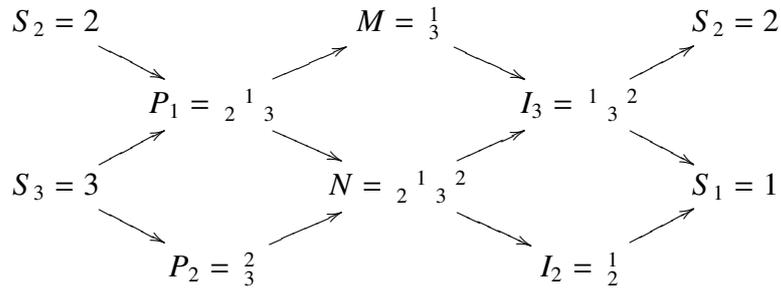
\xymatrix@=4mm{
S_2 = 2 \ar@{..}[rr]  \ar[dr] & &  M ={\begin{smallmatrix} 1 \\ 3 \end{smallmatrix}} \ar@{..}[rr]   \ar[dr]& & S_2 = 2 \\
& P_1 = {\begin{smallmatrix} & 1 &  \\ 2 & & 3 \end{smallmatrix}} \ar@{..}[rr] \ar[ur] \ar[dr]& & I_3 = {\begin{smallmatrix} 1&  &2  \\  & 3&  \end{smallmatrix}}\ar[ur] \ar[dr] & \\
S_3 = 3 \ar@{..}[rr] \ar[ur] \ar[dr] & & N = {\begin{smallmatrix} & 1 && 2  \\  2 && 3 & \end{smallmatrix}} \ar@{..}[rr] \ar[ur] \ar[dr] & & S_1 = 1\\
& P_2 = {\begin{smallmatrix} 2 \\ 3 \end{smallmatrix}} \ar@{..}[rr] \ar[ur] & & I_2 = {\begin{smallmatrix}  1   \\ 2  \end{smallmatrix}} \ar[ur]& 
}
$$
where the notation indicates which simple modules occur in the radical layers of the module, so $N = {\begin{smallmatrix} & 1 && 2  \\  2 && 3 & \end{smallmatrix}}$ is a module of length 4, of radical length 2, and with
$N/\rad N \simeq S_1 \amalg S_2$.

The following table gives, for each $\tau$-rigid indecomposable $U$,
a list of the indecomposable modules in $J(U)$ and an algebra $\Gamma_U$
such that $J(U)\simeq \module \Gamma$.

\setlength{\extrarowheight}{2.2mm}
\begin{center}
\begin{tabular}{|c|c|c| }
\hline 
$U$ & $J(U)$ & $\Gamma_U$ \\
\hline 
$S_2 = 2 $ & $\{ {\begin{smallmatrix} 2 \\ 3 \end{smallmatrix}},
{\begin{smallmatrix} 1&  &2  \\  & 3&  \end{smallmatrix}}, 1 \}$ & 
$k(\sbt \longrightarrow \sbt)$\\ \hline 
$S_3 = 3$ & $\{ 2, {\begin{smallmatrix} 1 \\ 2 \end{smallmatrix}}, 1 \}$ &
$k(\sbt  \longrightarrow \sbt )$ \\ \hline 
$P_1 = {\begin{smallmatrix} & 1 &  \\ 2 & & 3 \end{smallmatrix}}$ &
$\{ 3, {\begin{smallmatrix} 2 \\ 3 \end{smallmatrix}} , 2 \}$ & $k(\sbt  \longrightarrow \sbt )$ \\ \hline 
$P_2 = {\begin{smallmatrix} 2  \\  3 \end{smallmatrix}}$ &
$\{ 3, {\begin{smallmatrix} 1 \\ 3 \end{smallmatrix}} , 1 \}$ & $k(\sbt  \longrightarrow \sbt )$\\ \hline 
$M = {\begin{smallmatrix} 1  \\  3 \end{smallmatrix}}$ &
$\{ 3, {\begin{smallmatrix} & 1 & \\ 2 & & 3 \end{smallmatrix}} , {\begin{smallmatrix}  1  \\ 2  \end{smallmatrix}} \}$ & $k(\sbt  \longrightarrow \sbt )$\\ \hline 
$N = {\begin{smallmatrix} & 1 && 2  \\  2 && 3 & \end{smallmatrix}}$ &
$\{ {\begin{smallmatrix}  1  \\  3 \end{smallmatrix}} , {\begin{smallmatrix}  2  \\ 3  \end{smallmatrix}} \}$ & $k(\sbt \hspace{9mm} \sbt)$\\ \hline 
$I_2 = {\begin{smallmatrix} 1   \\ 2  \end{smallmatrix}}$ &
$\{ {\begin{smallmatrix} & 1 & \\ 2 & & 3 \end{smallmatrix}}, {\begin{smallmatrix} 1 \\ 3 \end{smallmatrix}}  , 
 {\begin{smallmatrix} & 1 && 2  \\  2 && 3 & \end{smallmatrix}},
 {\begin{smallmatrix} 1 & & 2  \\ & 3 &   \end{smallmatrix}}, 2 \}$ & 
 $k(\xymatrix{ \sbt \ar^{\alpha}@<2pt>[r] & \sbt \ar^{\beta}@<2pt>[l]}) / (\beta \alpha)$
\\ \hline 
$S_1$ &  $\{ {\begin{smallmatrix}  1  \\  2 \end{smallmatrix}} , {\begin{smallmatrix}  1  \\ 3  \end{smallmatrix}} \}$ & 
$k(\sbt \hspace{9mm} \sbt)$\\
\hline 
\end{tabular}
\end{center}

\bigskip

\noindent As remarked in the introduction, the module $I_3$ is {\em not} $\tau$-rigid, but as an object in the module category $J(S_2)$, it is $\tau$-rigid.
Hence $(P_2,I_3,S_2)$ is a $\tau$-exceptional sequence.
Since $I_3$ is projective in $J(S_2)$, we have that $(P_2,I_3[1],S_2)$ is also a $\tau$-exceptional sequence.

\noindent
We calculate the total number of signed $\tau$-exceptional sequences as follows.
For $$U \in \{S_2, S_3= P_3, P_2, P_1, M, P_3[1], P_2[1], P_1[1]\}$$ there are (by Example 1) 10 signed $\tau$-exceptional sequences of the form 
$(-,-,U)$. 
For $U \in \{N, S_1\}$ there are 8 signed $\tau$-exceptional sequences of the form 
$(-,-,U)$, while there are (by Example 2) 12 signed $\tau$-exceptional sequences of the form 
$(-,-,I_2)$. Hence, in total there are 108 signed $\tau$-exceptional sequences for
this algebra.

We conclude with examples illustrating how we compute
which signed $\tau$-exceptional sequence is
the image of a given support $\tau$-tilting object under our bijection.
Consider the ordered 
support $\tau$-rigid object $(M, I_2, P_1)$.
To compute $\E_{P_1} (M)$ we first note that $M$ is in $\Gen P_1$,
so $\E_{P_1} (M) = \rho(M)$.
Furthermore we have that $\PP_M$ is given by $P_2 \to P_1$ and so we
have the triangle (\ref{BCtrianglei}) in this case is:
$$P_2 \to P_1 \to \PP_M$$
and hence $\rho(M) = (f_{P_1}(P_2))[1] = P_2[1]$.
Similarly, we have that $\E_{P_1} (I_2) = \rho(I_2) = S_3[1]$.
The ordered 
support $\tau$-rigid object $(P_2[1], S_3[1])$ in $J(P_1)$ corresponds
according to the table of Example 1 to the signed $\tau$-exceptional sequence
$(S_2[1], S_3[1])$. Hence our bijection maps the 
ordered 
support $\tau$-rigid object $(M, I_2, P_1)$ to the 
signed $\tau$-exceptional sequence $(S_2[1], S_3[1], P_1)$.

Consider the ordered 
support $\tau$-rigid object $(M, P_1, I_2)$.
Note that $\Hom(I_2, M) = 0= \Hom(I_2, P_1)$, so that 
$\E_{I_2}(M) = M$ and $\E_{I_2}(P_1) = P_1$.
The ordered 
support $\tau$-rigid object $(M, P_1)$ in $J(I_2)$ corresponds,
according to the table of Example 2, to the signed $\tau$-exceptional sequence
$(S_2[1], P_1)$ in $J(I_2)$ (note that $S_2$ is projective in $J_{J(I_2)}(P_1)$). Hence our bijection maps the 
ordered 
support $\tau$-rigid object $(M, P_1, I_2)$ to the 
signed $\tau$-exceptional sequence $(S_2[1], P_1, I_2)$.

\end{document}